\documentclass[11pt]{article}

\def\eb{{\bf e}}

\def\h{{\bf h}}

\def\ubold{{\bf u}}
\def\vbold{{\bf v}}
\def\w{{\bf w}}
\def\x{{\bf x}}
\def\y{{\bf y}}
\def\z{{\bf z}}

\def\B{{\mathcal B}}
\def\C{{\mathcal C}}

\def\F{{\mathcal F}}

\def\M{{\mathcal M}}
\def\N{{\mathcal N}}

\def\S{{\mathcal S}}

\def\R{{\mathbb R}}

\def\al{\alpha}
\def\d{\delta}
\def\de{\delta}
\def\D{\Delta}
\def\e{\epsilon}

\def\ga{\gamma}

\def\l{\lambda}
\def\la{\lambda}

\def\om{\omega}
\def\OM{\Omega}

\def\s{\sigma}
\def\SI{\Sigma}

\def\Th{\Theta}

\def\balpha{{\boldsymbol \alpha}}

\def\bD{{\boldsymbol \D}}

\def\boldeta{{\boldsymbol \eta}}

\def\bth{{\boldsymbol \theta}}
\def\bphi{{\boldsymbol \phi}}

\def\bxi{{\boldsymbol \xi}}
\def\bzeta{{\boldsymbol \zeta}}

\def\tai{t \ap \infty}

\def\ap{\rightarrow}

\def\seq{\subseteq}

\def\bi{\{0,1\}}

\def\bp{\{-1,1\}}
\def\bz{{\bf 0}}

\def\fa{\; \forall}

\def\half{\frac{1}{2}}

\def\as{\mbox{ a.s.}}

\def\nm{\Vert}

\def\gJ{\nabla J}
\def\gJt{\gJ(\bth_t)}
\def\bxt{\bxi_{t+1}}
\def\bzt{\bzeta_{t+1}}
\def\bths{\bth^*}

\def\gJ{\nabla J}

\newcommand{\bc}{\begin{center}}
\newcommand{\ec}{\end{center}}
\newcommand{\be}{\begin{equation}}
\newcommand{\ee}{\end{equation}}
\newcommand{\bd}{\begin{displaymath}}
\newcommand{\ed}{\end{displaymath}}
\newcommand{\ba}{\begin{array}}
\newcommand{\ea}{\end{array}}
\newcommand{\ben}{\begin{enumerate}}
\newcommand{\een}{\end{enumerate}}
\newcommand{\bit}{\begin{itemize}}
\newcommand{\eit}{\end{itemize}}
\newcommand{\beq}{\begin{eqnarray}}
\newcommand{\eeq}{\end{eqnarray}}
\newcommand{\btab}{\begin{tabular}}
\newcommand{\etab}{\end{tabular}}
\newcommand{\bfig}{\begin{figure}}
\newcommand{\efig}{\end{figure}}
\newcommand{\btp}{\begin{tikzpicture}}
\newcommand{\etp}{\end{tikzpicture}}

\newcommand{\nmm}[1]{ \nm #1 \nm }
\newcommand{\nmeu}[1]{ \nm #1 \nm_2 }
\newcommand{\nmeusq}[1]{ \nm #1 \nm_2^2 }

\newcommand{\IP}[2]{ \langle #1 , #2 \rangle }

\def\nmsl1{\nm_{{\rm SL1}}}

\usepackage{graphicx}
  \DeclareGraphicsExtensions{.jpeg,.png,.jpg,.eps,.pdf}
\usepackage{tikz}
\usepackage{epstopdf}
\usepackage{hyperref}
\usepackage{amssymb,amsmath}

\def\bphH{\h_{t+1}}

\newtheorem{corollary}{Corollary}{\bf}{\it}
\newtheorem{definition}{Definition}{\bf}{\it}
{\bf}{\rm}
\newtheorem{lemma}{Lemma}{\bf}{\it}
\newtheorem{theorem}{Theorem}{\bf}{\it}
{\bf}{\it}
{\bf}{\it}
{\bf}{\rm}
\newtheorem{proof}{Proof}{\bf}{\rm}

\def\gJ{\nabla J}
\def\gJt{\gJ(\bth_t)}
\def\gJw{\gJ(\w_t)}

\def\xbb{\bar{\x}}
\def\bD{{\boldsymbol \D}}

\begin{document}

\title{
Convergence of the Stochastic Heavy Ball Method \\
With Approximate Gradients and/or Block Updating
}
\author{
Uday Kiran Reddy Tadipatri\thanks{University of Pennsylvania, ukreddy@seas.upenn.edu} 
\and 
M.\ Vidyasagar \thanks{Indian Institute of Technology Hyderabad, m.vidyasagar@iith.ac.in }}

\maketitle

\begin{abstract}

In this paper, we establish the convergence of the stochastic Heavy Ball
(SHB) algorithm under more general conditions than in the current literature.
Specifically,
(i) The stochastic gradient is permitted to be biased, and also, to
have conditional variance that grows over time (or iteration number).
This feature is essential when applying SHB with zeroth-order methods,
which use only two function evaluations to approximate the gradient.
In contrast, all existing papers assume that the stochastic gradient
is unbiased and/or has bounded conditional variance.
(ii) The step sizes are permitted to be random, which is essential when
applying SHB with block updating.
The sufficient conditions for convergence are stochastic analogs of the
well-known Robbins-Monro conditions.
This is in contrast to existing papers where more restrictive conditions
are imposed on the step size sequence.
(iii) Our analysis embraces not only convex functions, but also more general
functions that satisfy the PL (Polyak-{\L}ojasiewicz) and KL
(Kurdyka-{\L}ojasiewicz) conditions.
(iv) If the stochastic gradient is unbiased and has bounded variance,
and the objective function satisfies (PL), then the iterations of SHB match
the known best rates for convex functions.
(v) We establish the almost-sure convergence of the iterations,
as opposed to convergence in the mean or convergence in probability,
which is the case in much of the literature.
(vi) Each of the above convergence results continue to hold if 
full-coordinate updating is replaced by any one of three widely-used
updating methods.
In addition,
numerical computations are carried out to illustrate the above points.

\textbf{This paper is dedicated to the memory of Boris Teodorovich Polyak}
\end{abstract}

\section{Introduction}\label{sec:Intro}

In this paper, we study the well-known and widely-used ``Heavy Ball'' (HB)
method for convex and nonconvex optimization introduced by 
Polyak in \cite{Polyak-CMMP64}, and subsequently studied by several
researchers.
The objective of the present paper is to prove the \textit{almost sure}
convergence of this algorithm, when it is applied to a class of nonconvex
objective functions (which includes convex functions as a special case),
when the search direction is random -- the so-called \textbf{stochastic
gradient}.
The stochastic gradient is permitted to be biased, and/or to have a
conditional variance that grows without bound as a function of the
iteration counter.
The assumptions in this paper are the least restrictive in the literature.
In addition, we show that the iterations converge \textit{almost surely}.

\subsection{Literature Review}\label{ssec:12}

Since the literature on optimization is vast,
our review is limited to the papers that are directly
relevant to the specific class of algorithms studied here.
We first review ``momentum-based'' methods in general, and then
we review relevant papers in the literature that establish the
\textit{almost sure} convergence of various algorithms.

Suppose the objective function to be minimized is $\C^1$, and denote it by
$J(\cdot)$.
The general form of the SHB algorithm studied in 
\cite{Sebbouh-et-al-CoLT21} is as follows:
Choose an initial guess $\bth_0 \in \R^d$ (either in a deterministic or
a random fashion).
For $t \geq 0$, choose a random vector $\h_{t+1} \in \R^d$,
known as the \textbf{stochastic gradient},
and update $\bth_t$ using the formula
\be\label{eq:111}
\bth_{t+1} = \bth_t - \al_t \h_{t+1} + \mu_t (\bth_t - \bth_{t-1}) .
\ee
Here $\mu_t \in [0,1)$ is called the ``momentum coefficient,''
and $\al_t$ is the step size at time $t$.
If $\mu_t = 0$ for all $t$, then \eqref{eq:111} becomes the standard
Stochastic Gradient Descent (SGD) algorithm, which is studied here in
Section \ref{ssec:21}.
If
\be\label{eq:112}
\h_{t+1} = \gJt + \bxt ,
\ee
where $\bxt$ is a measurement error,
\eqref{eq:111} becomes ``stochastic gradient descent with
a momentum term.''
The usual assumption is that the error $\bxt$ has zero conditional mean
and bounded conditional variance.
However, in this paper we do permit more general stochastic gradients
than in \eqref{eq:112}.
The more general type of error is essential when the stochastic gradient
is determined using only function evaluations.
More details can be found in Section \ref{ssec:13}.

The Heavy Ball (HB)  method, introduced in \cite{Polyak-CMMP64}, is 
one of earliest ``momentum-based'' methods for optimization.
The phrase ``momentum-based'' is somewhat vague, but refers to
methods wherein the search direction at step $t$
depends not only on the current guess
$\bth_t$, but also on the previous guess $\bth_{t-1}$.
The algorithm introduced in \cite{Polyak-CMMP64} is
\be\label{eq:113}
\bth_{t+1} = \bth_t - \al \gJt + \mu (\bth_t - \bth_{t-1}) ,
\ee
which is of the form \eqref{eq:111} with $\h_{t+1} = \gJt$,
and $\al_t, \mu_t$ set equal to constants.
It is shown by Polyak that, if
$J(\bth)$ is quadratic of the form $(1/2) \bth^\top A \bth
+ \IP{\vbold}{\bth} + c$ for some positive definite matrix $A$,
vector $\vbold$ and constant $c$,
then the HB method requires $1/\sqrt{R}$ fewer iterations
compared to the gradient descent method, provided $\mu$ is
chosen as $(\sqrt{R}-1)/(\sqrt{R}+1)$, where $R$ denotes the condition
number of $A$.

A subsequent and widely-used momentum-based method is Nesterov's
Accelerated Gradient (NAG) method \cite{Nesterov-Dokl83}.
In \cite{Sutskever-et-al-ICML13}, NAG is reformulated in a manner
that brings out the similarities as well as the differences with HB.
Specifically, the NAG algorithm can be written as
\be\label{eq:121}
\vbold_{t+1}^N = \mu_t \vbold_t^N - \al_t \gJ[\bth_t + \mu_t \vbold^N_t ] ,
\vbold_0^N = \bz ,
\ee
\be\label{eq:122}
\bth_{t+1} = \bth_t + \vbold^N_{t+1} .
\ee
These two equations can be combined into the single equation
\be\label{eq:123}
\bth_{t+1} = \bth_t - \al_t \gJ[\bth_t + \mu_t( \bth_t - \bth_{t-1} ) ]
+ \mu_t( \bth_t - \bth_{t-1} ) .
\ee
This can be put in the HB formulation \eqref{eq:111} if we define
\be\label{eq:124}
\h_{t+1} = \gJ[\bth_t + \mu_t( \bth_t - \bth_{t-1} ) ] .
\ee
In other words, in NAG the gradient is computed \textit{after} the
momentum correction term $\mu_t (\bth_t - \bth_{t-1})$ is added to $\bth_t$.
The paper \cite{Sutskever-et-al-ICML13}
also shows that NAG can be deployed successfully
for training deep neural networks.
Also, it is shown in \cite[Section 2.2]{Nesterov18}
that when $J(\cdot)$ is a smooth convex
function with a Lipschitz-continuous gradient, NAG
converges to the minimum at the rate of $O(t^{-2})$.
Moreover, no gradient-based algorithm can achieve a faster rate.
More details can be found in \cite[Section 7]{Bottou-et-al-SIAM18}.

Another relevant reference is \cite{Bengio-et-al-ICASSP13},
in which the NAG algorithm is reformulated in an equivalent form, namely
\be\label{eq:125}
\vbold_{t+1}^B = \mu_t \vbold^B_t - \al_t \gJ(\Th_t) , 
\vbold_0^B = \bz ,
\ee
\beq
\Th_{t+1} & = & \Th_t + (1 + \mu_t) \vbold^B_{t+1} - \mu_{t-1} \vbold_t^B 
\nonumber \\
& = & \Th_t + \mu_t \mu_{t-1} \vbold_t^B + (1 + \mu_t) \al_t 
\gJ(\Th_t) \label{eq:126} .
\eeq
If started with the initial guess $\bth_0 = \bz$, the trajectory of
this algorithm matches that in \cite{Sutskever-et-al-ICML13}
(which is just a reformulation of NAG) both at the start and in the
final phase of local convergence to the solution.
But the formulation in \cite{Bengio-et-al-ICASSP13} is closer to Polyak's
HB compared to NAG, because the gradient $\gJ(\cdot)$ is computed
at the current guess $\Th_t$, and not a shifted version of it.

A fruitful approach to analyzing the behavior of NAG is to study
an associated second-order ODE on $\R^d$.
This is done in \cite{Su-Boyd-Candes16}, when
the step size $\al$ is held constant, while the momentum coefficient
$\mu_t$ varies with time.
It is shown that the ``optimal'' schedule for $\mu_t$ is
$\mu_t = (t+2)/(t+5)$.
In \cite{Aujol-et-al-SIAMO19}, the rate of convergence of the ODE
is analyzed further by imposing additional structure on $J(\cdot)$,
such as the {\L}ojasiewicz property (defined in Section \ref{ssec:20}
below).
It is shown that, under such assumptions, it is possible for classical
steepest descent method to outperform NAG.
The Heavy Ball algorithm can also be analyzed via its own associated ODE,
which too is a second-order ODE in $\R^d$.
This ODE is analyzed when $J(\cdot)$ satisfies either the 
Polyak-{\L}ojasiewicz property \cite{Apido-et-al-JGO22},
or the {\L}ojasiewicz property \cite{Aujol-et-al-MP23}.
%
In all of the above formulations, it is assumed that the ``stochastic
gradient'' equals the true gradient $\gJt$; thus these models do not allow
for measurement errors.

Now we come to more recent research on HB.
Much of the initial work studied the behavior of HB when $J(\cdot)$
is quadratic; however, later work expanded the scope to cover
strictly convex or convex functions.
In much of the literature, attention is focused in the \textit{convergence
in expectation} of various algorithms; sometimes \textit{convergence in
probability} is also studied.
However, any stochastic algorithm generates \textit{one sample path}
of a stochastic process.
Thus it is very useful to know that almost all sample paths converge
to the desired limit.
In the review paper \cite{Bottou-et-al-SIAM18},
the emphasis is almost exclusively on convergence in expectation.
SHB and SNAG are discussed in \cite[Section 7]{Bottou-et-al-SIAM18}.
Other research on the convergence of HB
(without establishing almost sure convergence)
is summarized very well on page 3 of \cite{Sebbouh-et-al-CoLT21}
and Section 1.1 of \cite{Liu-Yuan-arxiv22}.
However, for the convenience of the reader, we summarize some of the
relevant papers.

To proceed further, we introduce a little notation regarding the
stochastic gradient $\h_{t+1}$ introduced in \eqref{eq:111}.
Let $\bth_0^t$ denote $( \bth_0 , \cdots , \bth_t )$, and similarly
$\h_1^t := ( \h_1 , \cdots , \h_t )$.
(Note that there is no $\h_0$.)
To cater to the possibility of the step size $\al_t$
being random, we also define $\balpha_0^t$ as above.
This situation arises when we study ``block'' updating in the later part
of the paper.
Suppose that all of these are random variables on some underlying
probability space $(\OM , \SI , P)$.\footnote{The reader is referred to
\cite{Durrett19} for all concepts related to stochastic processes,
conditional expectations, etc.}
Let $\F_t$ denote the $\s$-algebra generated by
$\bth_0$, $\h_1^t$, $\balpha_0^t$.
Then it is clear that $\{ \F_t \}$ is a \textbf{filtration}, that is,
an increasing sequence of $\s$-algebras.
Moreover, if we denote the set of all random variables that are measurable
with respect to $\F_t$ by $\M(\F_t)$, then
$\bth_0^t , \h_1^t, \balpha_0^t \in \M(\F_t)$.
Observe that if the step size sequence is deterministic
(but possibly a function of $t$), then $\F_t$ is the $\s$-algebra generated by
$\bth_0$ and $\h_1^t$.

For future use, if $X \in \R^d$ is a random vector and $\{ \F_t \}$
is a filtration, we let $E_t(X)$ denote
$E( X | \F_t )$, the conditional expectation of $X$
with respect to $\F_t$.
Also, let $CV_t(X)$ denote the conditional variance of $X$, that is
\be\label{eq:127}
CV_t(X) := E_t( \nmeusq{ X - E_t(X) } )  = E_t( \nmeusq{X}) - \nmeusq{E_t(X)} .
\ee

For a stochastic gradient $\h_{t+1}$, define
\be\label{eq:128}
\z_t := E_t(\h_{t+1}) , \x_t := \z_t - \gJt , \bzt := \h_{t+1} - \z_t .
\ee
Thus $\x_t$ measures the ``bias'' of the stochastic gradient, that is,
the difference between the conditional expectation $E_t(\h_{t+1})$
and the true gradient $\gJt$.
From the ``tower'' property of conditional expectations
(see \cite{Williams91}), it follows that\footnote{Hereafter
we omit the phrase ``almost surely'' almost everywhere.}
\be\label{eq:129}
E_t(\bzt) = \bz \as, \fa t .
\ee
As a result,
\be\label{eq:1210}
CV_t(\h_{t+1}) = \nmeusq{\z_t} + E_t ( \nmeusq{\bzt} ) .
\ee
Note that $\x_t$ quantifies the difference between the
conditional expectation of the stochastic gradient, and the true gradient.
Thus, if $\x_t = \bz$, then
\be\label{eq:1210b}
E_t (\h_{t+1}) = \gJt ,
\ee
so that $\h_{t+1}$ is an unbiased estimate of the gradient.
In such a case, \eqref{eq:1210} simplifies to
\be\label{eq:1210c}
CV_t(\h_{t+1}) = E_t ( \nmeusq{\bzt} ) .
\ee
In much of the literature, the phrase ``stochastic gradient'' is used to
refer to the case where $\x_t = \bz$, and in addition,
there exists a finite constant $M$ such that
\be\label{eq:1210d}
CV_t(\h_{t+1}) = E_t ( \nmeusq{\bzt} ) \leq M^2 .
\ee
However, we will find it profitable to interpret the phrase more broadly,
and to introduce the definitions in \eqref{eq:128}.

Now we give a brief literature review in chronological order,
starting with papers that study the SGD, and then move to SHB.
Most papers only prove convergence in expectation.

In \cite{Nemirovski-et-al-SIAMO09}, the authors study the minimization of
convex functions that are not necessarily differentiable.
The assumption on the stochastic gradient is that $\z_t = E_t(\h_{t+1}) 
\in \partial J(\bth_t)$, where $\partial(\cdot)$ denotes the subgradient.
The paper analyzes the standard iteration in \eqref{eq:111} with
$\mu_t \equiv 0$, so that the algorithm under study is SGD
and not SHB.
The paper gives a very general analysis of the ``averaging'' approach
proposed by Polyak and Ruppert, and studied in \cite{Pol-Jud-SICOPT92}.
In \cite{Bach-Moul-NIPS11}, the authors extend some earlier results
from strictly convex functions to convex functions.
In \cite{Ghadimi-Lan-SIAMO13}, the authors study SGD where the gradient
makes use of ``zeroth-order stochastic function values,'' that is,
function values corrupted by noise.
These are used to compute approximate derivatives.
In this sense, the contents are in the spirit of
\cite{Kief-Wolf-AOMS52,Blum54}.
The paper \cite{Mertik-et-al-NeurIPS20} is one of the few to establish almost
sure convergence of the iterations, but under some very strong assumptions.
For example, $J(\cdot)$ is assumed to be $d$-times differentiable,
which can be a problem if $d$ is large.
The usual assumption elsewhere in the literature
is that $J \in \C^1$ no matter what $d$ is.
It is also assumed in \cite{Mertik-et-al-NeurIPS20}
that the gradient $\gJt$ is \textit{globally bounded},
which means that $J(\cdot)$ is restricted to \textit{linear} growth.

Now we come to papers that specifically study HB.
In \cite{Ghadimi-et-al-ECC15}, the authors analyze the HB algorithm
where $\h_{t+1} = \gJt$; thus there is no provision for measurement noise,
so that the algorithm being analyzed is HB and not SHB.
The function $J(\cdot)$ is assumed to be convex, and to have a globally
Lipschitz-continuous gradient.
The authors \textit{do not} show that $J(\bth_t)$ converges to the global
minimum of $J(\cdot)$.
Rather, they show that \textit{the average of the first $t$ iterations}
converges to the minimum value of the function.
In \cite{Gadat-et-al-EJS18}, the authors study the SHB for
some classes of nonconvex functions.
It is assumed that the stochastic gradient is unbiased, i.e., that
$E_t(\h_{t+1}) = \gJt$, so that $\x_t = \bz$ for all $t$.
The iterations are shown to converge to a minimum,
but at the cost of ``uniformly elliptic bounds'' on the measurement error
$\bzt$, which are very restrictive.
Finally, we mention \cite{Liu-et-al-NC23}, in which a ``unified''
algorithm is presented, which includes woth SHB and SNAG as special cases.
In that paper, only convergence in expectation is proved, and that too,
under the assumption that the stochastic gradient $\h_{t+1}$ is
unbiased and has finite conditional variance.

Now we discuss three papers that are most closely related to the present paper,
namely \cite{Ber-Tsi-SIAM00,Sebbouh-et-al-CoLT21,Liu-Yuan-arxiv22}.

Perhaps the closest in spirit to the present paper is the old paper 
\cite{Ber-Tsi-SIAM00}.
In that paper, the stochastic gradient $\h_{t+1}$ is allowed to be biased,
and the assumptions on the conditional variance of $\h_{t+1}$ are
similar to ours.
The authors also prove almost sure convergence.
The theorems in \cite{Ber-Tsi-SIAM00} do not apply to SHB, just to SGD;
moreover, it is unclear how their arguments can be adapted to handle SHB.
Nevertheless, it is an important paper.
Most of the complexity of the proofs in \cite{Ber-Tsi-SIAM00} arises
because the authors permit $J(\cdot)$
to be unbounded from below, whereas in most of the literature (including this
paper), it is assumed that $J^* := \inf_{\bth} J(\bth) > - \infty$.

In \cite{Sebbouh-et-al-CoLT21}, the objective function is an expected
value, of the form (\cite[Eq.\ (1)]{Sebbouh-et-al-CoLT21})
\bd
J(\bth) = E_{\w \sim P} F(\bth,\w) .
\ed
The function $F(\cdot,\w)$ is convex for each $\w$, and its gradient
is Lipschitz-continuous with constant $L_\w \leq L$ for all $\w$.
Thus $J(\cdot)$ is also convex, and $\gJ(\cdot)$ is also $L$-Lipschitz
continuous.
The stochastic gradient is chosen as (\cite[Eq.\ (SHB)]{Sebbouh-et-al-CoLT21})
\bd
\h_{t+1} = \nabla_{\bth_t} F(\w_{t+1},\bth_t) ,
\ed
where $\w_{t+1}$ is chosen i.i.d.\ with distribution $P$.
Effectively this means that in \eqref{eq:128},
$\z_t = E_t(\h_{t+1}) = \gJt$, so that $\x_t = \bz$.
In other words, the stochastic gradient is \textit{unbiased}.
Also, it is assumed that, for some constant $\s^2$, we have that
(\cite[Eq.\ (5)]{Sebbouh-et-al-CoLT21})
\be\label{eq:1210a}
CV_t(\bzt) \leq 4L (J(\bth_t) - J^* ) + \s^2 ,
\ee
where $J^*$ is the infimum of $J(\cdot)$.
Thus the hypotheses are more restrictive than \eqref{eq:222} and \eqref{eq:223},
which are the assumptions in the present paper.

In \cite{Sebbouh-et-al-CoLT21}, it is suggested how to convert \eqref{eq:111}
above to \textit{two} equations which do not contain any ``delayed'' terms.
Specifically, the authors iteratively define
\be\label{eq:412}
\la_{t+1} = \frac{\la_t}{\mu_t} - 1 ,
\eta_t = (1 + \la_{t+1}) \al_t 
\ee
In the above, in principle the quantity $\la_0$ is not specified and 
can be chosen by the user.
If we now define
\be\label{eq:413}
\w_{t+1} = \w_t - \eta_t \h_{t+1} ,
\ee
\be\label{eq:414}
\bth_{t+1} = \frac{ \la_{t+1} }{ 1 + \la_{t+1} } \bth_t 
+ \frac{1}{1 + \la_{t+1} } \w_{t+1} ,
\ee
then $\bth_{t+1}$ satisfies \eqref{eq:111}.
Note that one can write \eqref{eq:414} as
\bd
\bth_{t+1} = \bth_t + \frac{1}{1 + \la_{t+1} } ( \w_{t+1} - \bth_t )
= \bth_t + \frac{1}{1 + \la_{t+1} } (\w_t - \bth_t - \eta_t \h_{t+1} ) ,
\ed
or equivalently as
\be\label{eq:415}
\bth_{t+1} 
= \bth_t + \frac{1}{1 + \la_{t+1} } (\w_t - \bth_t ) - \al_t \h_{t+1}  .
\ee 
Then the equations \eqref{eq:413} and \eqref{eq:415} together
resemble an SGD in the joint variable $(\bth_t,\w_t)$.
Thus in principle the standard results on the convergence of the SGD
can be used to analyze \eqref{eq:413}--\eqref{eq:414}.

In \cite{Sebbouh-et-al-CoLT21}, the authors choose $\la_0 = 0$, and
\be\label{eq:415b}
\la_t = \frac{S_{t-1}}{\eta_t} , \mbox{ where } 
S_t = \sum_{\tau = 0}^t \eta_\tau .
\ee
In Condition 1, it is assumed that
\be\label{eq:1211}
\eta_{t+1} \leq \eta_t ,
\sum_{t=0}^\infty \eta_t^2 \s^2 < \infty ,
\sum_{t=0}^\infty \eta_t = \infty,
\sum_{t=0}^\infty ( \eta_t /S_t)  = \infty,
\ee
Under these assumptions, it is shown in Theorem 13 that
\bd
J(\bth_t) - J^* = o(1/S_{t-1}) .
\ed
Note that the conditions in \eqref{eq:1211} are more restrictive
than the standard Robbins-Monro conditions in terms of the amended
step size sequence $\{ \eta_t \}$.
First, the step size sequence is assumed to be decreasing, and second,
since the sequence $\{ S_t \}$ is strictly increasing, the last condition
in \eqref{eq:1211} is more restrictive than the divergence of the sum
of $\eta_t$.

In addition, the main challenge in the above approach is that there is no
obvious and verifiable relationship between the original parameters $\al_t$
(step size) and $\mu_t$ (the momentum parameter), and the
convergence conditions \eqref{eq:1211}.
In particular, even if the
original step size sequence $\{ \al_t \}$ satisfies the Robbins-Monro
conditions
\be\label{eq:1212}
\sum_{t=0}^\infty \al_t^2 < \infty ,
\sum_{t=0}^\infty \al_t = \infty,
\ee
the sequence $\{ 1 + \la_{t+1} \}$
might increase too rapidly for the sequence $\{ \eta_t \}$ to satisfy
\eqref{eq:1211}.
This is why the authors \textit{begin} with the sequence $\{ \eta_t \}$.
Clearly, it would be desirable to state the convergence conditions
directly in terms of the step size sequence $\{ \al_t \}$ and the
momentum sequence $\{ \mu_t \}$.
In the present paper, we show that when the momentum parameter is a constant,
then the standard Robbins-Monro conditions \eqref{eq:1212} are sufficient for 
convergence.

In \cite{Liu-Yuan-arxiv22} the authors study the case where the
objective function is either strongly convex, or nonconvex 
with a Lipschitz-continuous gradient.
Unlike in \cite{Sebbouh-et-al-CoLT21}, these authors assume that
$\mu_t$ is constant.
In this case, $\la_t$  is a constant, and $\eta_t$ is a constant
multiple of $\al_t$.
Hence
\bd
\la = \frac{\mu}{1-\mu} , 1 + \la = \frac{1}{1 - \mu} , 
\frac{\la}{1+\la} = \mu ,
\eta_t = \frac{\al_t}{1 - \mu} .
\ed
Therefore \eqref{eq:1211} become just the standard Robbins-Monro conditions
on the step size $\al_t$.

In this case, \eqref{eq:413} and \eqref{eq:415} become
\be\label{eq:415a}
\bth_{t+1} = \mu \bth_t + (1 - \mu) \w_t - \al_t \h_{t+1} , \quad
\w_{t+1} = \w_t - \frac{\al_t}{1 - \mu} \h_{t+1} .
\ee
In \cite{Liu-Yuan-arxiv22}, the authors do not use the above equations.
Instead, they define
\bd
\vbold_t = \bth_t - \bth_{t-1}, \y_t = \bth_t + \frac{\mu}{1-\mu} \vbold_t ,
\ed
and show that these two quantities satisfy the recursions
\be\label{eq:416}
\vbold_{t+1} = \mu \vbold_t - \al_t \h_{t+1} , 
\y_{t+1} = \y_t - \frac{\al_t}{1-\mu} \h_{t+1} .
\ee
It is assumed that $\h_{t+1}$ is unbiased (i.e., that $\x_t = \bz$),
and that the variance satisfies the ``Expected Smoothness''
condition proposed in \cite{Khaled-Rich-arxiv20},
As shown in \cite{MV-RLK-SGD-JOTA24}, the expected smoothness
condition is more restrictive than the conditions assumed here.
When the objective function is strictly convex, 
the authors study only the case where $\al_t = \Th(t^{1-\phi})$ for some
$\phi \in (0,1/2)$, and show that
\be\label{eq:1213}
J(\bth) - J^* = o(1/(t^{1-\e})), \fa \e \in (2 \phi, 1) .
\ee
We build upon this this approach in the present paper.
When the function is not strongly convex, and just has a Lipschitz-continuous
gradient, it is assumed that 
\be\label{eq:1214}
\al_{t+1} \leq \al_t , \quad
\sum_{t=0}^\infty (\al_t/S_{t-1}) = \infty,
\mbox{ where } S_t = \sum_{\tau = 0}^t \al_t .
\ee
Thus \eqref{eq:1214} is basically the same as \eqref{eq:1211} when
$\mu_t$ is a constant.
In this case the authors prove a weaker conclusion than \eqref{eq:1213}, namely
\be\label{eq:1215}
\min_{0 \leq \tau \leq t} \nmeusq{\gJ(\bth_\tau)} = o(1/S_{t-1} ) .
\ee

In all of the papers discussed until now, \textit{every component}
of $\bth_t$ is updated at each step $t$, according to \eqref{eq:111}.
This might be referred to as ``synchronous'' updating,
though this terminology is not very standard.
At the other end of the spectrum lies ``coordinate'' updating,
in which the update \eqref{eq:111} is applied to \textit{exactly one}
randomly chosen component of $\bth_t$.
Note that the phrase ``coordinate updating'' is not very standard.
However, the phrase ``coordinate gradient descent'' is quite standard.
When the measurements are noise-free, the behavior of coordinate gradient
descent is analyzed \cite{Wright15} for convex functions, and in
\cite{Bach-et-al-aisats19} for a class of nonconvex functions.
However, the results in these papers do not apply when the gradient
measurement is corrupted by noise.
In contrast, the results presented here can cope with noisy measurements.
In-between synchronous updating and coordinate updating lies what
we choose to call ``block updating'' (or Block Coordinate Descent (BCD)),
in which, at each step $t$,
some subset $S_t \seq [d]$ is chosen, and only those components of
$\bth_{t,i}, i \in S_t$ are updated using \eqref{eq:111}.
Observe that in block updating, both the cardinality of the set $S_t$
as well as the elements can be random.
The convergence of block updating with error-free measurements
has been studied in \cite{Nesterov-SICOPT12,
Richtarik-Takac-MP12, Richtarik-Takac-MP15}.
In \cite{Lu-Xiao-arxiv13}, the authors 
provide a probabilistic convergence result, based
on the Nesterov's framework \cite{Nesterov-SICOPT12}.
The study was limited to smooth convex functions and bounded noise variance.
Convergence of block updating in SGD for nonconvex functions has not
been studied much.
In \cite{Xu-Yin-SICOPT15}, the convergence of BCD is proved for
nonconvex functions under bounded noise in the measurements.
Block updating has gained a lot
of attention in  distributed ML \cite{Niu-et-al-arxiv11}, broadly categorized
into two main algorithms: Synchronous SGD (updates are performed
one after another node) \cite{Chen-et-al-iclr16} and Asynchronous SGD (ASGD)
(random updates by any node at any time) \cite{Niu-et-al-arxiv11,
Xie-et-al-ICML20}.
In short, the present paper addresses a combination of
issues that are not found in the existing literature, to the best
of the authors' knowledge.

\subsection{Contributions of the Paper}\label{ssec:13}

For technical reasons explained below, we restrict attention to the case where
the momentum coefficient $\mu_t$ in \eqref{eq:111} is constant, as in
\cite{Liu-Yuan-arxiv22}, but unlike \cite{Sebbouh-et-al-CoLT21}.
However, unlike both these papers, we permit the step size $\al_t$
to be \textit{random}.
This is crucial for studying block-updating, as described in Section
\ref{sec:Meta} below.

Now we discuss the contributions of the present paper.
\bit
\item All of our convergence results hold for \textit{arbitrary
and possibly random} step size
sequences $\{ \al_t \}$ that satisfy stochastic analogs of
the Robbins-Monro conditions \cite{Robbins-Monro51} or
the Kiefer-Wolfowitz-Blum conditions \cite{Kief-Wolf-AOMS52,Blum54}.
In contrast, in \cite{Sebbouh-et-al-CoLT21,Liu-Yuan-arxiv22},
the step sizes are deterministic and need to satisfy
more stringent assumptions as in \eqref{eq:1211} above.
In both these papers, it is possible to choose the step size as
$\al_t = 1/(t+1)^s$ for a suitable exponent $s$.
However, it may be advantageous to have a convergence proof that
requires nothing more than the standard Robbins-Monro conditions.
\item Our assumptions on the stochastic gradient are the less restrictive
than those in the current literature.
Specifically, we permit the stochastic gradient to be biased, and
also permit the bias to grow linearly with respect to $\gJt$.
Similarly, we permit the conditional variance of the stochastic gradient
to increase with respect to $t$, and also at a rate of $J(\bth_t)$.
In contrast, in both \cite{Sebbouh-et-al-CoLT21,Liu-Yuan-arxiv22},
the stochastic gradient is assumed to be unbiased.
Our assumption is weaker than \cite[Eq.\ (5)]{Sebbouh-et-al-CoLT21},
which implies that the conditional variance of the stochastic gradient
is bounded both with respect to $t$ as well as $\bth_t$.
Also, our assumption is weaker than the ``Expected Smoothness''
condition proposed in \cite{Khaled-Rich-arxiv20}, and is assumed in
\cite{Liu-Yuan-arxiv22}.
The Expected Smoothness assumption is 
the weakest assumption in the literature to date, prior to our paper.
\item As a result of these relaxed assumptions, the theory presented here
can be used to establish the convergence of the
SHB algorithm when the stochastic gradient $\h_{t+1}$
is computed using only function valuations
(sometimes referred to as a ``gradient-free'' or ``zeroth-order'' method).
In particular, we show that the SPSA (Simultaneous Perturbation
Stochastic Approximation) introduced in 
\cite{Spall-TAES98,Sadegh-Spall-TAC98,Hernandez-Spall-ACC14} works also
when a momentum term is introduced.
So far as we are aware, this is a first.
\item We establish the \textit{almost sure convergence} of the algorithm
when a stochastic gradient is used instead of the true gradient.
While there is some literature on \textit{convergence in expectation},
there are not many results on almost sure convergence.
\item We study the minimization of a class of \textit{nonconvex} objective
functions, which is more general than those studied thus far.
Specifically, when the objective function satisfies an analog of the
Kurdyka-{\L}ojasiewicz property, we establish the almost sure convergence
of the objective function to its minimum value.
Under the stronger Polyak-{\L}ojasiewicz property, we not only establish
almost sure convergence, but also bounds on the \textit{rate} of
convergence.
\item
We study \textbf{Block Stochastic Gradient Descent (BSGD)} and
\textbf{Block Stochastic Heavy Ball} algorithms, where in,
at each iteration, \textit{some but not necessarily all} components
of the current guess are updated.
We prove a ``meta-theorem'' to the effect that, when an SGD algorithm
converges with full coordinate update by virtue of satisfying the
sufficient conditions in \cite{MV-RLK-SGD-JOTA24},
the same algorithm continues to converge with block updating as well.
We prove the convergence of the SHB algorithm by converting it to
an SGD algorithm with more variables, as in \cite{Liu-Yuan-arxiv22}. 
Consequently, the convergence of Block SHB also follows readily.
\eit

\subsection{Organization of the Paper}\label{ssec:14}

The paper is organized as follows:
In Section \ref{sec:SGD}, we reprise some relevant results from
\cite{MV-RLK-SGD-arxiv23,MV-RLK-SGD-JOTA24} on the convergence
of the Stochastic Gradient Descent (SGD) algorithm.
As it turns out, the problem formulation put forward in these papers provides
a framework that also embraces the Stochastic Heavy Ball (SHB) algorithm,
either with full-coordinate or block-updating.
In Section \ref{sec:SHB}, we state precisely the version of SHB that is
under study here, and then proceed to prove our main results with
\textit{full} coordinate updating.
In Section \ref{sec:Meta}, we state and prove a ``meta'' theorem
for the convergence of block-updating in general.
While the meta-theorem is applied here to SHB alone, the meta-theorem
is quite useful by itself, in our views.
In Section \ref{sec:Num}, we present numerical results on the application
of not just SHB but a variety of algorithms, on three distinct objective
functions, out of which two are not convex.
Finally, in Section \ref{sec:Conc}, we summarize our contributions,
and mention some research topics that merit further investigation.

\section{Reprise of Relevant Results for SGD}\label{sec:SGD}

In this section, we restate some relevant results from \cite{MV-RLK-SGD-JOTA24}
on the convergence of the Stochastic Gradient Descent (SGD) algorithm
\be\label{eq:130}
\bth_{t+1} = \bth_t - \al_t \h_{t+1} ,
\ee
where $\h_{t+1}$ is the stochastic gradient and $\al_t$ is the step size.
The proofs of the cited results can be found in that reference.
These results form the basis for the convergence results for the SHB
algorithm in Section \ref{sec:SHB}, and also for
the ``meta'' theorem in Section \ref{sec:Meta}
on the convergence of SGD or SHB with block updating.

\subsection{Standing Assumptions and Their Significance}\label{ssec:20}

The theory in \cite{MV-RLK-SGD-JOTA24}
applies to a class of smooth nonconvex functions, as well as to
all smooth convex functions.
Moreover, the error models used there are the least restrictive among those
found in the literature to date.
To make these points precise, we begin by discussing the class(es) of
functions under study, and then the error models.

We begin with two ``standing'' assumptions on the objective function
$J(\cdot)$.
These assumptions are standard in the literature, and
assumed to hold in the remainder of the paper.
\ben
\item[(S1)] $J(\cdot)$ is $\C^1$, and $\gJ(\cdot)$ is globally 
Lipschitz-continuous with constant $L$.
\item[(S2)] $J(\cdot)$ is bounded below, and the infimum is attained.
Thus
\bd
J^* := \inf_{\bth \in \R^d} J(\bth)
\ed
is well-defined, and $J^* > -\infty$.
Moreover, the set
\be\label{eq:131}
S_J := \{ \bth : J(\bth) = J^* \}
\ee
is nonempty.
\textit{By redefining $J(\cdot)$ if necessary, hereafter it is assumed
that $J^* = 0$.}
\een

Before proceeding further, we draw the reader's attention to the
following useful result.

\begin{lemma}\label{lemma:11}
Suppose (S1) holds, and that $J^*  > - \infty$.
Then
\be\label{eq:134}
\nmeusq{ \nabla J(\bth) } \leq 2L [ J(\bth) - J^* ].
\ee
\end{lemma}

This result is Lemma 4.1 of \cite{MV-RLK-SGD-JOTA24}.
For future use, this bound is referred to as the Gradient Growth (GG)
property.

For functions that satisfy (S1) and (S2), we delineate various properties.
Note that different theorems assume different properties on $J(\cdot)$,
which in turn lead to different conclusions.
Define, as usual,
\bd
\rho(\bth) := \inf_{\bphi \in S_J} \nmeu{\bth - \bphi} 
\ed
to be the distance from $\bth$ to the set of minimizers $S_J$.
Also, we define a \textbf{function of Class $\B$} to be a map $\psi: \R_+ \ap
\R_+$ such that $\psi(0) = 0$, and 
\bd
\inf_{\e \leq x \leq M} \psi(x) > 0
\ed
whenever $0 < \e \leq M < \infty$.
With the aid of this definition, we introduce two function classes.
\ben
\item[(PL)] There exists a constant $K$ such that
\be\label{eq:132}
\nmeusq{\gJ(\bth)} \geq K J(\bth) , \fa \bth \in \R^d .
\ee
\item[(KL')] There exists a function $\psi(\cdot)$ of Class $\B$
such that
\be\label{eq:133}
\nmeu{\gJ(\bth)} \geq \psi(J(\bth) , \fa \bth \in \R^d .
\ee
\een
Finally, we introduce one last property.
\ben
\item[(NSC)]
This property consists of the following assumptions, taken together.
\bit
\item The function $J(\cdot)$ attains its infimum.
Therefore the set $S_J$ defined in \eqref{eq:131} is nonempty.
\item
The function $J(\cdot)$ has compact level sets.
Thus for every constant $c \in (0,\infty)$, the level set
\bd
L_J(c) := \{ \bth \in \R^d : J(\bth) \leq c \}
\ed
is compact.
\item
There exists a function $\eta(\cdot)$ of Class $\B$ such that
\be\label{eq:135}
\rho(\bth) \leq \eta(J(\bth)) , \fa \bth \in \R^d .
\ee
\eit
\een

Next we discuss the significance of these assumptions, as well as the
nomenclature.

PL stands for the Polyak-{\L}ojasiewicz condition.
In \cite{Polyak-UCMMP63}, Polyak introduced \eqref{eq:132},
and showed that it is sufficient to ensure that iterations converge at
a ``linear'' (or geometric) rate to a global minimum, whether or not
$J(\cdot)$ is convex.
Note that \eqref{eq:132} can also be rewritten as
\bd
\nmeu{\gJ(\bth)} \geq K^{1/2} [J(\bth)]^{1/2} , \fa \bth \in \R^d .
\ed
In \cite{Loja63}, {\L}ojasiewicz introduced a more general condition
\be\label{eq:136}
\nmeu{J(\bth)} \geq C [ J(\bth)]^r , \fa \bth \in \R^d ,
\ee
for some constant $C$ and some exponent $r \in [1/2,1)$.
Note that in the present paper, we use only the Polyak condition
\eqref{eq:132}.

In \cite{Kurdyka98}, Kurdyka proposed a more general inequality than
\eqref{eq:136}, 
namely: There exist a constant $c > 0$ and a function
$v: [0,c) \ap \R$ which is $\C^1$ on $(0,c)$, such that $v'(x) > 0$
for all $x \in (0,c)$, and
\be\label{eq:138}
\nmeu{\gJ(\bth)} \geq [ v'(J(\bth) ]^{-1} .
\ee
In particular, if $v(x) = x^{1-r}$ for some $r \in (0,1)$, then
\eqref{eq:138} becomes \eqref{eq:136} with $C = 1/(1-r)$.
For this reason, \eqref{eq:138} is sometimes referred to as the
Kurdyka-{\L}ojasiewicz (KL) inequality.
See for example \cite{BDLM-TAMS10}.
In our case, we don't require the right side to be a differentiable
function; rather we require only that
it be a function of Class $\B$ of $J(\bth)$.
Hence we choose to call this condition as (KL'), to suggest that it is
similar to, but weaker than, the KL condition.
Note that, under (PL) or (KL'), $\gJ(\bth) = \bz$ implies that
$J(\bth) = 0$, i.e., that every stationary point is also a global minimum.
Thus any function that satisfies either (PL) or (KL') is ``invex''
as defined in \cite{Hanson-JMAA81}.
See \cite{Karimi-et-al16} for an excellent survey of these topics.

Finally, (NSC) stands for ``Near Strong Convexity.''
A function $J(\cdot)$ is said to be \textbf{$R$-strongly convex} if
\bd
J(\bth) \geq J(\bphi) + \IP{\gJ(\bphi)}{\bth - \bphi}
+ \frac{R}{2} \nmeusq{\bth - \bphi} , \fa \bth , \bphi \in \R^d .
\ed
Note that an $R$-strongly convex function has a unique global
minimizer $\bths$.
If we relax the assumption of strong convexity and
ask only that the above relation holds for all $\bphi \in S_J$
(which need not be a singleton set),
then, after noting that $\gJ(\bphi) = \bz$ for all $\bphi \in \S_J$,
the above bound becomes
\be\label{eq:138a}
J(\bth) \geq \frac{R}{2} \nmeusq{\bth - \bphi} , \fa \bth \in \R^d ,
\fa \bphi \in S_J .
\ee
Since
\bd
\rho(\bth) = \inf_{\bphi \in S_J} \nmeu{\bth - \bphi} ,
\ed
It follows that if $J(\cdot)$ satisfies \eqref{eq:138a}, then
\bd
[ (2/R) J(\bth)]^{1/2} \geq \rho(\bth) , \fa \bth \in \R^d .
\ed
In property (NSC), the left side of the above is changed
to $\eta(J(\bth))$ for some function $\eta(\cdot)$ of Class $\B$.
Thus it is a very mild assumption.
A consequence of the (NSC) property is that,
whenever $J(\cdot)$ satisfies (NSC), and $J(\bth_t) \ap 0$
as $\tai$, we can conclude that $\rho(\bth_t) \ap 0$, i.e., that
$\bth_t$ approaches the set of minimizers of $J(\cdot)$.

Now we reprise a very useful bound, namely \cite[Eq.\ (2.4)]{Ber-Tsi-SIAM00}.
This is a workhorse of several proofs in this paper.

\begin{lemma}\label{lemma:12}
Suppose $J: \R^d$ is $\C^1$, and that $\gJ(\cdot)$ is $L$-Lipschitz
continuous.
Then
\be\label{eq:139}
J(\bth+\bphi) \leq J(\bth) + \IP{\gJ(\bth)}{\bphi} + \frac{L}{2} 
\nmeusq{\bphi} , \fa \bth , \bphi \in \R^d .
\ee
\end{lemma}

\subsection{Relevant Results on the Convergence of SGD}\label{ssec:21}

In this subsection we quote, without proof, some relevant results
from \cite{MV-RLK-SGD-JOTA24} on the convergence of the Stochastic
Gradient Descent (SGD) algorithm.
The reader can consult this source for the proofs of the cited results.
The contents of this subsection are relevant because we prove the
convergence of SHB, either with full or with Block updating,
by invoking the results presented here.

A fundamental result in the convergence of stochastic processes
is the classic ``almost supermartingale'' theorem due to
Robbins and Siegmund \cite[Theorem 1]{Robb-Sieg71}.
It is also found in \cite{BMP90} and in \cite{Fran-Gram21}.
The Robbins-Siegmund theorem states the following:

\begin{lemma}\label{lemma:21}
Suppose $\{ z_t \} , \{ f_t \} , \{ g_t \} , \{ h_t \}$ are
stochastic processes taking values in $[0,\infty)$, adapted to some
filtration $\{ \F_t \}$, satisfying
\be\label{eq:211}
E_t( z_{t+1} ) \leq (1 + f_t) z_t + g_t - h_t \as, \fa t ,
\ee
where, as before, $E_t(z_{t+1})$ is a shorthand for $E(z_{t+1} | \F_t )$.
Then, on the set
\bd
\OM_0 := \{ \om \in \OM : \sum_{t=0}^\infty f_t(\om) < \infty \}
\cap \{ \om : \sum_{t=0}^\infty g_t(\om) < \infty \} ,
\ed
we have that $\lim_{\tai} z_t$ exists, and in addition,
$\sum_{t=0}^\infty h_t(\om) < \infty$.
In particular, if $P(\OM_0) = 1$, then $\{ z_t \}$ is bounded
almost surely, and $\sum_{t=0}^\infty h_t(\om) < \infty$ almost surely.
\end{lemma}

The following theorem is a straight-forward, but useful extension
of Lemma \ref{lemma:21}.
It is Theorem 5.1 of \cite{MV-RLK-SGD-JOTA24},
and can be used to establish
the convergence of stochastic gradient methods for nonconvex functions.

\begin{theorem}\label{thm:21}
Suppose $\{ z_t \} , \{ f_t \} , \{ g_t \} , \{ h_t \}, \{ \al_t \}$ are
$[0,\infty)$-valued stochastic processes
defined on some probability space $(\OM,\SI,P)$, and
adapted to some filtration $\{ \F_t \}$.
Suppose further that
\be\label{eq:212}
E_t(z_{t+1} ) \leq (1 + f_t) z_t + g_t - \al_t h_t \as, \fa t .
\ee
Define
\be\label{eq:213}
\OM_0 := \{ \om \in \OM : \sum_{t=0}^\infty f_t(\om) < \infty \mbox{ and }
\sum_{t=0}^\infty g_t(\om) < \infty \} ,
\ee
\be\label{eq:213f}
\OM_1 := \{ \sum_{t=0}^\infty \al_t(\om) = \infty \} .
\ee
Then
\ben
\item Suppose that $P(\OM_0) = 1$.
Then the sequence $\{ z_t \}$ is bounded almost surely, and
there exists a random variable $W$ defined on $(\OM,\SI,P)$ such that
$z_t(\om) \ap W(\om)$ almost surely.
\item
Suppose that, in addition to $P(\OM_0) = 1$, it is also true that
$P(\OM_1) = 1$.
Then
\be\label{eq:214}
\liminf_{\tai} h_t(\om) = 0  \fa \om \in \OM_0 \cap \OM_1 .
\ee
Further, suppose there exists a function $\eta(\cdot)$ of Class $\B$
such that $h_t(\om) \geq \eta(z_t(\om))$ for all $\om \in \OM_0$.
Then $z_t(\om) \ap 0$ as $\tai$ for all $\om \in \OM_0 \cap \OM_1$.
\een
\end{theorem}

Theorem \ref{thm:21} allows us to infer convergence, but does not
provide any information about the \textit{rate} of convergence.
Now we define the concept of a rate of convergence of stochastic
processes, following a similar definition in \cite{Liu-Yuan-arxiv22}.

\begin{definition}\label{def:order}
Suppose $\{ Y_t \}$ is a stochastic process, and $\{ f_t \}$
is a sequence of positive numbers.
We say that
\ben
\item $Y_t = O(f_t)$ if $\{ Y_t / f_t \}$ is bounded almost surely.
\item $Y_t = \OM(f_t)$ if $Y_t$ is positive almost surely, and
$\{ f_t/Y_t \}$ is bounded almost surely.
\item $Y_t = \Th(f_t)$ if $Y_t$ is both $O(f_t)$ and $\OM(f_t)$.
\item $Y_t = o(f_t)$ if $Y_t /f_t \ap 0$ almost surely as $\tai$.
\een
\end{definition}

With this definition, the following theorem holds; it is Theorem 5.2 of
\cite{MV-RLK-SGD-JOTA24}.
Similar results can be found in \cite{Liu-Yuan-arxiv22}.

\begin{theorem}\label{thm:22}
Suppose $\{ z_t \} , \{ f_t \} , \{ g_t \} , \{ \al_t \}$ are
stochastic processes defined on some probability space $(\OM,\SI,P)$,
taking values in $[0,\infty)$, adapted to some
filtration $\{ \F_t \}$.
Suppose further that
\be\label{eq:213a}
E_t(z_{t+1} ) \leq (1 + f_t) z_t + g_t - \al_t z_t \fa t ,
\ee
where
\bd
\sum_{t=0}^\infty f_t(\om) < \infty ,
\sum_{t=0}^\infty g_t(\om) < \infty ,
\sum_{t=0}^\infty \al_t(\om) = \infty .
\ed
Then $z_t = o(t^{-\l})$ for every $\l \in (0,1]$ such that
there exists a finite $T > 0$ such that
\be\label{eq:213b}
\al_t(\om) - \l t^{-1} \geq 0 \fa t \geq T ,
\ee
and in addition
\be\label{eq:213c}
\sum_{t=0}^\infty (t+1)^\la g_t(\om) < \infty ,
\sum_{t=0}^\infty [ \al_t(\om) - \la t^{-1} ] = \infty .
\ee
\end{theorem}

Next, we consider the SGD (that is,
\eqref{eq:111} with $\mu_t = 0$ for all $t$).
Thus
\be\label{eq:221}
\bth_{t+1} = \bth_t - \al_t \h_{t+1} ,
\ee
where $\h_{t+1}$ is the stochastic gradient.
Recall the various quantities defined in \eqref{eq:128}.
Suppose the stochastic gradient satisfies the following assumptions:
There exist sequences of constants $\{ B_t \}$ and $\{ M_t \}$ such that
\be\label{eq:222}
\nmeu{\x_t} \leq B_t [ 1 + \nmeu{\gJt} ] , \fa \bth_t \in \R^d , \fa t ,
\ee
\be\label{eq:223}
E_t ( \nmeusq{\bzt} ) = CV_t(\h_{t+1}) \leq M_t^2 [ 1 + J(\bth_t) ] ,
\fa \bth_t \in \R^d , \fa t .
\ee
With these assumptions, we can state the following result, which is
Theorem 6.1 of \cite{MV-RLK-SGD-JOTA24}.

\begin{theorem}\label{thm:23}
Suppose the objective function $J(\cdot)$ satisfies the standing assumptions
(S1) and (S2), as well as Property (GG).
Suppose further that the stochastic gradient $\h_{t+1}$ satisfies
\eqref{eq:222} and \eqref{eq:223}.
With these assumptions, we have the following conclusions:\footnote{All
hypotheses and conclusions hold almost surely.}
\ben
\item Suppose that
\be\label{eq:224}
\sum_{t=0}^\infty \al_t^2 < \infty , 
\sum_{t=0}^\infty \al_t B_t  < \infty ,
\sum_{t=0}^\infty \al_t^2 M_t^2 < \infty .
\ee
Then $\{ \gJ(\bth_t) \}$ and $\{ J(\bth_t) \}$ are bounded, and in addition,
$J(\bth_t)$ converges to some random variable as $\tai$.
\item If in addition $J(\cdot)$ satisfies (KL'), and
\be\label{eq:225}
\sum_{t=0}^\infty \al_t = \infty ,
\ee 
then $J(\bth) \ap 0$ and $\gJ(\bth_t) \ap \bz$ as $\tai$.
\item Suppose that in addition to (KL'), $J(\cdot)$ also satisfies (NSC),
and that \eqref{eq:222} and \eqref{eq:223} both hold.
Then $\rho(\bth_t) \ap 0$ as $\tai$.
\een
\end{theorem}

Finally, by strengthening the hypothesis from Property (KL') to Property (PL),
we can state a bound on the \textit{rate of convergence} of SGD.
It is Theorem 6.2 of \cite{MV-RLK-SGD-JOTA24}.

\begin{theorem}\label{thm:24}
Let various symbols be as in Theorem \ref{thm:23}.
Suppose $J(\cdot)$ satisfies the standing assumptions (S1) and (S2),
as well as (GG).
Suppose in addition that $J(\cdot)$ satisfies property (PL),
and that \eqref{eq:224} and \eqref{eq:225} hold.
Further, suppose there exist constants $\ga > 0$ and $\de \geq 0$ such
that the constants $B_t$ and $M_t$ in \eqref{eq:222} and \eqref{eq:223}
satisfy\footnote{Since $t^{-\ga}$ is undefined when $t = 0$,
we really mean $(t+1)^{-\ga}$.
The same applies elsewhere also.}
\bd
B_t = O(t^{-\ga}), M_t = O(t^\de) ,
\ed
where we take $\ga = 1$ if $B_t = 0$ for all sufficiently large $t$,
and $\de = 0$ if $M_t$ is bounded.
Choose the step-size sequence $\{ \al_t \}$ as
$O(t^{-(1-\phi)})$ and $\OM(t^{-(1-C)})$
where $\phi$ and $C$ are chosen to satisfy
\bd
0 < \phi < \min \{ 0.5 - \de , \ga \} , C \in (0,\phi] .
\ed
Define
\be\label{eq:226}
\nu := \min \{ 1 - 2( \phi + \de) , \ga - \phi \} .
\ee
Then $\nmeusq{\gJt} = o(t^{-\la})$ and $J(\bth_t) = o(t^{-\la})$
for every $\la \in (0,\nu)$.
In particular, by choosing $\phi$ very small, it follows that
$\nmeusq{\gJt} = o(t^{-\la})$ and $J(\bth_t) = o(t^{-\la})$ whenever
\be\label{eq:227}
\la < \min \{ 1 - 2 \de , \ga \} .
\ee
\end{theorem}

\section{Convergence Theorems for the SHB Algorithm}\label{sec:SHB}

\subsection{Preliminaries}\label{ssec:41}

In this section we state and prove a convergence theorem for the Stochastic
Heavy Ball (SHB) algorithm with full coordinate update.
This is achieved by formulating SHB as an instance of SGD in an enlarged
variable space.

\subsection{Convergence of the SHB Algorithm Under Full Coordinate
Updating}\label{ssec:42}

In this subsection, we state our main results on the convergence of the
Stochastic Heavy Ball (SHB) algorithm under full coordinate updating.
Since several versions of the algorithm are studied in the literature
review, we now state explicitly the specific version being studied here.
The algorithm is
\be\label{eq:421}
\bth_{t+1} = \bth_t + \mu (\bth_t - \bth_{t-1}) - \al_t \h_{t+1} .
\ee
To analyze the algorithm, we rewrite \eqref{eq:421} equivalently as
in \eqref{eq:415a}, again reprised for the convenience of the reader.
\be\label{eq:424}
\bth_{t+1} = \mu \bth_t + (1 - \mu) \w_t - \al_t \h_{t+1} , \quad
\w_{t+1} = \w_t - \frac{\al_t}{1 - \mu} \h_{t+1} .
\ee

The following assumptions are made
\bit
\item The momentum term $\mu$ is \textit{constant}, and satisfies
$\mu \in [0,1)$.
\item The updating formula \eqref{eq:421} is applied to \textit{every
coordinate} of $\bth_t$ at each time $t+1$.
Thus we study full coordinate update.
The case of block updating is taken up in Section \ref{sec:Meta}.
\item The step size $\al_t$ is permitted to be random, and belongs
almost surely to $(0,\infty)$.
\eit

Let $\F_t$ be the $\s$-algebra generated by $\bth_0$, $\h_1^t$, and
if applicable, $\al_0^t$.
Let $E_t(X)$ denote the conditional expectation $E(X | \F_t )$.
In order to prove convergence, it is assumed that the stochastic gradient
$\h_{t+1}$ satisfies \eqref{eq:222} and \eqref{eq:223},
reprised here for the convenience of the reader.
There exist sequences of constants $\{ B_t \}$ and $\{ M_t \}$ such that
\be\label{eq:422}
\nmeu{\x_t} \leq B_t [ 1 + \nmeu{\gJt} ] , \fa \bth_t \in \R^d , \fa t ,
\ee
\be\label{eq:423}
E_t ( \nmeusq{\bzt} ) = CV_t(\h_{t+1}) \leq M_t^2 [ 1 + J(\bth_t) ] ,
\fa \bth_t \in \R^d , \fa t .
\ee

Now we state the convergence theorems.
The first theorem assures convergence when $J(\cdot)$ satisfies
the (KL') property, whereas the second theorem contains bounds
on the rate of convergence when $J(\cdot)$ satisfies the stronger
(PL) property.
Note that all hypotheses and conclusions hold almost surely.

\begin{theorem}\label{thm:421}
Suppose $J(\cdot)$ satisfies the assumptions (S1) and (S2),
while $\h_{t+1}$ satisfies \eqref{eq:422} and \eqref{eq:423}.
Then we have the following conclusions.
\ben
\item Suppose
\be\label{eq:425}   
\sum_{t=0}^\infty \al_t^2 < \infty ,
\sum_{t=0}^\infty \al_t B_t  < \infty ,
\sum_{t=0}^\infty \al_t^2 M_t^2 < \infty .
\ee
Then $\{ \gJ(\bth_t) \}$ and $\{ J(\bth_t) \}$ are bounded, and in addition,
$J(\bth_t)$ converges to some random variable as $\tai$.
\item If in addition $J(\cdot)$ satisfies (KL'), and
\be\label{eq:426}
\sum_{t=0}^\infty \al_t = \infty ,
\ee 
then $J(\bth) \ap 0$ and $\gJ(\bth_t) \ap \bz$ as $\tai$.      
\item Suppose that both \eqref{eq:425} and \eqref{eq:426} hold, and
in addition to (KL'), $J(\cdot)$ also satisfies (NSC).
Then $\rho(\bth_t) \ap 0$ as $\tai$.
\een
\end{theorem}

\begin{proof}
The first step is to ``decouple'' \eqref{eq:424} by a simple linear
transformation of the variables.
Define a new variable
\bd
\ubold_t := \bth_t - \w_t .
\ed
Then \eqref{eq:424} leads to
\beq
\ubold_{t+1} & = & \bth_{t+1} - \w_{t+1}
= \mu \bth_t - \mu \w_t + \frac{\mu}{1-\mu} \al_t \h_{t+1} \nonumber \\
& = & \mu \ubold_t +  \frac{\mu}{1-\mu} \h_{t+1} , \label{eq:433}
\eeq
while $\w_{t+1}$ continues to be governed by the second equation in
\eqref{eq:424}.
Note the ``sign inversion'' of the coefficient of $\h_{t+1}$.
It can be eliminated by changing the definition of $\ubold_t$ to $\w_t - \bth_t$,
which leads to slightly more messy formulas.
As we shall see, this ``sign inversion'' does not affect anything.

Because the proof is very elaborate, we state the general philosophy
up-front to guide the reader.
We use the ``Lyapunov function'' $J(\w_t) + \nmeusq{\ubold_t}$,
and then derive a bound in the form
\be\label{eq:433a}
\begin{split}
E_t( J(\w_{t+1}) + \nmeusq{\ubold_{t+1}} ) & \leq
J(\w_t) + \nmeusq{\ubold_t} + R_t - F_t \\
&- \left( 1 - \frac{\mu^2}{2} \right) \nmeusq{\ubold_t} - \al_t \nmeusq{\gJ(\w_t)} ,
\end{split}
\ee
where
\be\label{eq:433b}
R_t = f_t( J(\w_t) + \nmeusq{\ubold_t} ) + g_t ,
\ee
with $\{ f_t \} , \{ g_t \}$ being summable sequences, and
$F_t$ is a quadratic form in $\nmeu{\gJ(\w_t)}$ and $\nmeu{\ubold_t}$
which is positive definite for sufficiently large $t$.
Suppose $T$ is chosen such that $F_t \geq 0$ for all $t \geq T$
(and note that $T$ could be path-dependent).
Then we can start the analysis of \eqref{eq:433a} from time $T$,
and drop the term $F_t$ thereafter.
Then we can apply the Robbins-Siegmund theorem (Lemma \ref{lemma:21})
to \eqref{eq:433a} without the $F_t$ term,
which leads to the desired conclusions.

The first step is to convert the bounds \eqref{eq:422} and \eqref{eq:423},
which are stated in terms of $\bth_t$, to terms involving $\w_t$ and $\ubold_t$.
For this purpose, we use the $L$-Lipschitz continuity of $\gJ(\cdot)$.
Hence
\be\label{eq:437}
\nmeu{\gJt - \gJw} \leq L \nmeu{\bth_t - \w_t} = L \nmeu{\ubold_t} .
\ee
Define
\bd
\xbb_t = \z_t - \gJw = E_t(\h_{t+1}) - \gJw .
\ed
Then
\be\label{eq:437a}
\begin{split}
\nmeu{\xbb_t} & \leq \nmeu{ \z_t - \gJt } + \nmeu{\gJt - \gJw} \\
& \leq B_t( 1 + \nmeu{\gJt} ) + L \nmeu{\ubold_t} \\
& \leq B_t( 1 + \nmeu{\gJw} ) + L \nmeu{\ubold_t} ) + L \nmeu{\ubold_t} .
\end{split}
\ee
Now we use the obvious inequality $x \leq (1+x^2)/2$,
and invoke Lemma \ref{lemma:11} which states that $\nmeusq{\gJw}
\leq 2L J(\w_t)$.
Hence we can rewrite \eqref{eq:437a} as
\be\label{eq:438}
\begin{split}
\nmeu{\xbb_t} &\leq L \nmeu{\ubold_t}
+ B_t [ 1.5 + 0.5 \nmeusq{\gJw} ) + L \nmeu{\ubold_t} ] \\
&\leq L \nmeu{\ubold_t} + B_t [ 1 + L J(\w_t) + L \nmeu{\ubold_t} ] .
\end{split}
\ee
Note that we have replaced $1.5$ and $0.5$ by $2$ and $1$ in the above,
to avoid dealing with fractions.

Next we recast the bound on $E_t(\nmeusq{\bzt})$ in terms of $\w_t$
and $\ubold_t$.
For this purpose we use \cite[Eq.\ (2.4)]{Ber-Tsi-SIAM00}, 
stated here as Lemma \ref{lemma:12}, which gives
\be\label{eq:439}
J(\bth_t) = J( \w_t + \ubold_t ) \leq J(\w_t) + \IP{\gJw}{\ubold_t}
+ \frac{L}{2} \nmeusq{\ubold_t} .
\ee
Combining the above with \eqref{eq:423}, we get
\be\label{eq:4310}
\begin{split}
E_t( \nmeusq{\bzt}) & \leq M_t^2 ( 1 + J(\bth_t) ) \\
& \leq M_t^2 [ 1 + J(\w_t) + \IP{\gJw}{\ubold_t}
+ \frac{L}{2} \nmeusq{\ubold_t} ] .
\end{split}
\ee
By Schwarz' inequality and Lemma \ref{lemma:11}, we have that
\begin{eqnarray*}
\IP{\gJw}{\ubold_t} & \leq & \nmeu{\gJw} \cdot \nmeu{\ubold_t} \\
& \leq & \half [ \nmeusq{\gJw} + \nmeusq{\ubold_t} ] \\
& \leq & \half [ 2L J(\w_t) + \nmeusq{\ubold_t} ] .
\end{eqnarray*}
Substituting into \eqref{eq:4310} gives
\be\label{eq:4311}
E_t( \nmeusq{\bzt} ) 
\leq  M_t^2 (1+L) [ J(\w_t) + \half \nmeusq{\ubold_t} ] + M_t^2 .
\ee
We also generate an upper bound for $\nmeusq{\z_t}$ for later use, starting
with \eqref{eq:438}, and using Lemma \ref{lemma:11}.
\be\label{eq:4312}
\begin{split}
\nmeusq{\z_t} &= \nmeusq{ \gJw + \xbb_t } 
= \nmeusq{\gJw} + 2 \IP{\gJw}{\xbb_t} + \nmeusq{\xbb_t} \\
& \leq 2 \nmeusq{\gJw} + 2 \nmeusq{\xbb_t} 
\leq 4L J(\w_t) + 2 \nmeusq{\xbb_t} .
\end{split}
\ee
Now, from \eqref{eq:437a}, it follows that
\be\label{eq:4313}
\begin{split}
\nmeusq{\xbb_t} & \leq B_t^2 ( 1 + \nmeu{\gJw} )^2
+ L^2 ( 1+B_t^2) \nmeusq{\ubold_t}  \\
&+ 2 B_t(1+B_t) L \nmeu{\gJw} \cdot \nmeu{\ubold_t} .
\end{split}
\ee
By combining all of these bounds into the expression
\bd
E_t( \nmeusq{\h_{t+1}} ) = \nmeusq{\z_t} + E_t ( \nmeusq{\bzt} ),
\ed
we can obtain for the left side.
We will not however write it out for the time being.

After all these preliminary steps, we come to the key steps, namely,
to find upper bounds for $E_t( J(\w_{t+1}) )$ and $E_t( \nmeusq{\ubold_{t+1}})$.
First, 
\begin{eqnarray*}
\nmeusq{\ubold_{t+1}} & = & \left\nm \mu \ubold_t + \frac{\mu}{1-\mu} \al_t \h_{t+1} 
\right\nm_2^2 \\
& = & \mu^2 \nmeusq{\ubold_t} + \frac{2 \mu^2 }{1-\mu} \al_t \IP{\ubold_t}{\h_{t+1}}
+ \frac{\mu^2}{(1-\mu)^2} \al_t^2 \nmeusq{\h_{t+1}} .
\end{eqnarray*}
Since $\h_{t+1} = \gJw + \xbb_t + \bzt$, it follows that
\begin{eqnarray*}
E_t( \nmeusq{\ubold_{t+1}} ) & = & \mu^2 \nmeusq{\ubold_t} 
+ \frac{2 \mu^2 }{1-\mu} \al_t \IP{\ubold_t}{\gJw} \\
& + & \frac{2 \mu^2}{1-\mu} \al_t \IP{\ubold_t}{\xbb_t} 
+ \frac{\mu^2}{(1-\mu)^2} \al_t^2 [ \nmeusq{\z_t} + E_t(\nmeusq{\bzt}) ] .
\end{eqnarray*}
Using Schwarz' inequality and \eqref{eq:437a} gives
\be\label{eq:4314}
E_t( \nmeusq{\ubold_{t+1}} ) \leq \mu^2 \nmeusq{\ubold_t}
+ \frac{2 \mu^2 }{1-\mu} \al_t \nmeu{\gJw} \cdot \nmeu{\ubold_t} \\
+ \frac{2 \mu^2}{(1-\mu)^2} \al_t^2 \nmeusq{\ubold_t} + R_{1,t} ,
\ee
where $R_{1,t}$ consists of a bound for all the remaining terms, obtained
from \eqref{eq:437a}, \eqref{eq:438}, \eqref{eq:4311} and \eqref{eq:4312}.
Specifically,
\begin{eqnarray*}
R_{1,t} & = & \frac{2 \mu^2}{1-\mu} \al_t B_t [ 2 + L J(\w_t) + L \nmeu{\ubold_t} ]
+ \frac{\mu^2}{(1-\mu)^2} \al_t^2 [2 J(\w_t) + 2 \nmeusq{\xbb_t} ] \\
& + & \frac{\mu^2}{(1-\mu)^2} \al_t^2 M_t^2 
[ (1+L) ( J(\w_t) + \nmeusq{\ubold_t} ) + 1 ] .
\end{eqnarray*}
By assumption, \eqref{eq:425} holds.
In particular, the summability of $\al_t^2$ implies that $\al_t \ap 0$
as $\tai$, and as a result that $\al_t$ is bounded.
Combined with the summability of $\al_t B_t$, this shows that $\al_t^2 B_t$
is also summable.
Finally, since every absolutely summable sequence is also square-summable,
the summability of $\al_t B_t$ implies that $\al_t^2 B_t^2$ is also summable.
Note that the expression above for $R_{1,t}$
involves only the summable sequences
$\{ \al_t^2 \}$, $\{ \al_t B_t \}$ and $\{ \al_t^2 M_t^2 \}$.
Hence one can find a bound
\bd
R_{1,t} \leq f_{1,t} ( J(\w_t) + \nmeusq{\ubold_t} ) + g_{1,t} ,
\ed
where both $\{ f_{1,t} \}$ and $\{ g_{1,t} \}$ are summable.

A bound for $E_t(J(\w_{t+1}))$ can be derived using an entirely similar 
approach.
From Lemma \ref{lemma:12}, we have that
\begin{eqnarray*}
J(\w_{t+1}) & = & J\left( \w_t - \frac{\al_t}{1-\mu} \h_{t+1} \right) \\
& \leq & J(\w_t) - \frac{\al_t}{1-\mu} \IP{\gJw}{\h_{t+1}}
+ \frac{L}{2(1-\mu)^2} \al_t^2 \nmeusq{\h_{t+1}} .
\end{eqnarray*}
Therefore
\be\label{eq:4315}
\begin{split}
E_t(J(\w_{t+1})) & \leq J(\w_t) - \frac{\al_t}{1-\mu} \nmeusq{\gJw} 
- \frac{\al_t}{1-\mu} \IP{\gJw}{\xbb_t} \\
&+ \frac{L}{2(1-\mu)^2} \al_t^2 [ \nmeusq{\z_t} + E_t(\nmeusq{\bzt}) ] .
\end{split}
\ee
Now by applying \eqref{eq:438}, we get
\be\label{eq:4316}
\left| \frac{\al_t}{1-\mu} \IP{\gJw}{\xbb_t} \right|
\leq \frac{L}{1-\mu} \al_t \nmeu{\gJw} \cdot \nmeu{\ubold_t} + H_t ,
\ee
where $H_t$ equals $\al_t B_t$ multipliying some terms, and $\{ \al_t B_t \}$
is summable.
Hence
\begin{eqnarray*}
H_t + \frac{L}{2(1-\mu)^2} \al_t^2 \nmeusq{\h_{t+1}} & \leq & R_{2,t} \\
& \leq & f_{2,t} ( J(\w_t) + \nmeusq{\ubold_t} ) + g_{2,t} ,
\end{eqnarray*}
where the sequences $\{ f_{2,t} \} , \{ g_{2,t} \}$ are summable.
Adding \eqref{eq:4314} and \eqref{eq:4315} gives
\be\label{eq:4317}
\begin{split}
E_t( J(\w_{t+1}) + \nmeusq{\ubold_{t+1}} ) & \leq J(\w_t) + \nmeusq{\ubold_t} \\
& - \nmeusq{\ubold_t^2} \left[ 1 - \mu^2 + \frac{2 \mu^2 L}{1-\mu} \al_t \right] 
+ \nmeu{\gJw} \cdot \nmeu{\ubold_t} \frac{2 \mu^2 + L}{1-\mu} \al_t \\
& - \frac{\al_t} {1-\mu} \nmeusq{ \gJw} + R_t ,
\end{split}
\ee
where $R_t = R_{1,t} + R_{2,t}$ is of the form \eqref{eq:433b}.
To proceed further, let us define a quadratic form
\be\label{eq:4318}
\begin{split}
F_t &= \nmeusq{\ubold_t} \left[ \frac{1-\mu^2}{2} - \frac{2 \mu^2 L}{1-\mu} \al_t 
\right] + \nmeu{\gJw} \cdot \nmeu{\ubold_t} \left[ \frac{2 \mu^2 + L}{1-\mu} 
\right] \al_t \\
&+ \nmeusq{\gJw} \frac{\mu}{1-\mu} \al_t .
\end{split}
\ee
If we now split the term $((1 - \mu^2)/2) \nmeusq{\ubold_t}$ into two
equal parts, we can rewrite \eqref{eq:4317} as
\be\label{eq:4320}
\begin{split}
E_t( J(\w_{t+1}) + \nmeusq{\ubold_{t+1}} ) & \leq J(\w_t) + \nmeusq{\ubold_t} \\
&- \left( \frac{1-\mu^2}{2} \right) \nmeusq{\ubold_t} - \al_t \nmeusq{\gJw}
- F_t + R_t .
\end{split}
\ee
It is now shown that $F_t$ is a positive definite form for $t$ sufficiently
large; specifically, there exists a $T < \infty$ such that
$F_t \geq 0$ for all $t \geq $T.
Suppose we succeed in proving this.
Since we can always start our analysis of \eqref{eq:4317} starting at
time $T$, we can write
\be\label{eq:4321}
\begin{split}
E_t( J(\w_{t+1}) + \nmeusq{\ubold_{t+1}} ) & \leq J(\w_t) + \nmeusq{\ubold_t} \\
&- \left( \frac{1-\mu^2}{2} \right) \nmeusq{\ubold_t} - \al_t \nmeusq{\gJw}
+ R_t , \fa t \geq T .
\end{split}  
\ee
In other words, the term $-F_t$ is gone.
Now \eqref{eq:4321} is in a form to which the Robbins-Siegmund theorem
(Lemma \ref{lemma:11}) can be applied.
So let us now establish the positive definiteness of the quadratic form
for sufficiently large $t$.
Note that
\bd
F_t = \left[ \ba{cc} \nmeu{\ubold_t} \\ \nmeu{\gJw} \ea \right]^\top K_t
\left[ \ba{cc} \nmeu{\ubold_t} \\ \nmeu{\gJw} \ea \right] ,
\ed
where
\bd
K_t = \left[ \ba{cc} \frac{1-\mu^2}{2} - \frac{2 \mu^2 L}{1-\mu} \al_t & 
\frac{2 \mu^2 + L }{2(1-\mu)} \al_t \\
\frac{2 \mu^2 + L }{2(1-\mu)} \al_t & \frac{\mu}{1-\mu} \al_t \ea \right] .
\ed
Note that $K_t$ is of the form
\bd
K_t = \left[ \ba{cc} \frac{1-\mu^2}{2} - a \al_t & b \al_t \\ b \al_t & d \al_t
\ea \right] 
\ed
for suitable positive
constants $a,b,d$ which need not be written out explicitly.
A symmetric $2 \times 2$ matrix is positive definite if (and only if)
its trace and its determinant are both positive.
In this case
\bd
\mathrm{tr}(K_t) = \frac{1-\mu^2}{2} - (a-d) \al_t , \quad
\mathrm{det}(K_t) = \frac{1-\mu^2}{2} d \al_t - ( ad+bc) \al_t^2 .
\ed
Since, by hypothesis, $\sum_{t=0}^\infty \al_t^2 < \infty$,
it follows that $\al_t \ap 0$ as $\tai$.
Therefore $\mathrm{tr}(K_t) > 0$ for sufficiently large $t$.
Moreover, $\al_t^2$ approaches zero faster than $\al_t$.
This in turn implies that $\mathrm{det}(K_t) > 0$ for sufficiently large $t$.
Hence we conclude that $K_t$ is a positive definite matrix for sufficiently
large $t$.

Since it has already been established that $R_t$ has the form \eqref{eq:433b}
where the sequences $\{ f_t \}$, $\{ g_t \}$ are summable, we can now
apply Theorem \ref{thm:21} to \eqref{eq:4321}.

We begin wih Item 1.
Note that all statements hold ``almost surely,'' so this qualifier is not
repeated each time.
Suppose \eqref{eq:425} holds.
Then the following conclusions follow from Theorem \ref{thm:21}:
\bit
\item $J(\w_t) + \nmeusq{\ubold_t}$ is bounded.
Moreover, there is a random variable $X$ such that
$J(\w_t) + \nmeusq{\ubold_t} \ap X$ (almost surely) as $\tai$.
\item Further, almost surely
\be\label{eq:4322}
\sum_{t=0}^\infty \left( \frac{1-\mu^2}{2} \right) 
\nmeusq{\ubold_t} + \al_t \nmeusq{\gJw} < \infty .
\ee
\eit
Since the summands in \eqref{eq:4322} are both nonnegative, and $(1 - \mu^2)/2$
is just a constant, it follows that
\be\label{eq:4323}
\sum_{t=0}^\infty \nmeusq{\ubold_t} < \infty ,
\ee
\be\label{eq:4324}
\sum_{t=0}^\infty \al_t \nmeusq{\gJw} < \infty .
\ee
Now \eqref{eq:4323} implies that $\nmeusq{\ubold_t} \ap 0$ as $\tai$,
i.e., that $\ubold_t \ap \bz$ as $\tai$.
In turn, if $J(\w_t) + \nmeusq{\ubold_t} \ap X$, then $J(\w_t) \ap X$ as $\tai$.

Now recall that $\bth_t = \w_t + \ubold_t$.
Since $J(\cdot)$ is continuous and $\ubold_t \ap \bz$,
it follows that $J(\bth_t) \ap X$ as $\tai$.
The boundedness of $\{ J(\bth_t) \}$ follows from it being a convergent
sequence.
Finally, the boundedness of $\{ \gJt \}$ follows from Lemma \ref{lemma:11}.
Thus we have established Item 1.

Next we address Item 2 of the theorem.
The hypotheses are that, in addition to \eqref{eq:425}, \eqref{eq:426}
also holds, and $J(\cdot)$ satisfies Property (KL').
Then by definition there exists
a function $\psi : \R \ap \R$ in Class $\B$ such that
$\nmeu{\gJt} \geq \psi(J(\bth_t))$.
Recall that all the stochastic processes are defined on some underlying
probability space $(\OM,\SI,P)$.
Define
\bd
\OM_0 := \{ \om \in \OM : J(\bth(\om)) \ap X(\om) 
\& \nmeusq{\ubold_t(\om)} \ap 0 \} ,
\ed
\bd
\OM_1 := \{ \om \in \OM : \sum_{t=0}^\infty \al_t(\om) = \infty \} .
\ed
Note that if the step sizes are deterministic, then $\OM_1 = \OM$.
Define $\OM_2 = \OM_0 \cap \OM_1$, and note that $P(\OM_2) = 1$,
by Item 1.

The objective is to show that $X(\om) = 0$ for all $\om \in \OM_2$.
Once this is done, it would follow from Lemma \ref{lemma:11} that
\bd
\nmeu{\gJ(\bth_t(\om))} \leq [2L J (\bth_t(\om))]^{1/2}
\ap 0 \mbox{ as } \tai , \fa \om \in \OM_2 .
\ed
Accordingly, suppose that, for some $\om \in \OM_0$, we have that $X(\om) > 0$,
say $X(\om) = 2\e$.
Define
\bd
G(\om) := \sup_t J(\bth_t(\om)) .
\ed
Then $G(\om) < \infty$ because $\{ J(\bth_t(\om)) \}$ is a convergent sequence.
Define
\bd
\d := \half \inf_{\e \leq r \leq G(\om)} \psi(r) .
\ed
Then $\d > 0$ because $\psi(\cdot)$ is a function of Class $\B$.
Now choose a $T_0 < \infty$ such that $J(\bth(\om)) \geq \e$
for all $t \geq T_0$.
By the (KL') property, it follows that
\bd
\nmeu{\gJ(\bth(\om))} \geq 2 \d , \fa t \geq T_0 .
\ed
Next, choose $T_1 < \infty$ such that $\nmeu{\ubold_t(\om)} \leq \d/L$
for all $t \geq T_1$, and define $T_2 = \min \{ T_0, T_1 \}$.
Then it follows from the Lipschitz continuity of $\gJ(\cdot)$ that
\be\label{eq:4325}
\nmeu{\gJ(\w_t(\om))} \geq \nmeu{\gJ(\bth_t(\om))} - L \nmeu{\ubold_t(\om)}
\geq \d, \fa t \geq T_2 .
\ee
On the other hand, because $\om \in \OM_2$, we have that
\be\label{eq:4326}
\sum_{t=T_2} \al_t(\om) = \infty .
\ee
Thus \eqref{eq:4325} and \eqref{eq:4326} together imply that
\bd
\sum_{t=T_2}^\infty \al_t \nmeusq{ \gJt } = \infty .
\ed
Since this contradicts \eqref{eq:4324}, we conclude that no such
$\om \in \OM_2$ can exist.
In other words $X(\om) = 0$ for all $\om \in \OM_2$.
This establishes Item 2.

Item 3 is a ready consequence of Item 2 and Property (NSC).
If $\{ J(\bth_t) \}$ is bounded, then the fact that $J(\cdot)$
has compact level sets means that $\{ \bth_t \}$ is bounded.
Then the fact that $J(\bth_t) \ap 0$ as $\tai$ implies that $\rho(\bth_t) \ap 0$
as $\tai$; in other words, the distance from the iterate $\bth_t$ to
the set $S_J$ of global minima approaches zero.
Note that it is \textit{not} assumed that $S_J$ consists of a singleton.
\end{proof}

\begin{theorem}\label{thm:423}
Let various symbols be as in Theorem \ref{thm:421}.
Suppose $J(\cdot)$ satisfies the standing assumptions (S1) and (S2)
and also property (PL),
and that \eqref{eq:425} and \eqref{eq:426} hold.
Further, suppose there exist constants $\ga > 0$ and $\de \geq 0$ such
that
\bd
\mu_t = O(t^{-\ga}), \quad
M_t = O(t^\de) , \fa t \geq 1 ,
\ed
where we take $\ga = 1$ if $\mu_t = 0$ for all sufficiently large $t$,
and $\de = 0$ if $M_t$ is bounded.
Choose the step-size sequence $\{ \al_t \}$ as
$O(t^{-(1-\phi)})$ and $\OM(t^{-(1-C)})$
where $\phi$ and $C$ are chosen to satisfy
\be\label{eq:4331a}
0 < \phi < \min \{ 0.5 - \de , \ga \} , \quad
C \in (0,\phi] .
\ee
Define
\be\label{eq:4332}
\nu := \min \{ 1 - 2( \phi + \de) , \ga - \phi \} .
\ee
Then $\nmeusq{\gJt} = o(t^{-\la})$ and $J(\bth_t) = o(t^{-\la})$
for every $\la \in (0,\nu)$.
In particular, by choosing $\phi$ very small, it follows that
$\nmeusq{\gJt} = o(t^{-\la})$ and $J(\bth_t) = o(t^{-\la})$ whenever
\be\label{eq:4333}
\la < \min \{ 1 - 2 \de , \ga \} .
\ee
\end{theorem}

\begin{proof}
The proof, based on Theorem \ref{thm:421}, is basically the same as that
of Theorem 6.2 of \cite{MV-RLK-SGD-arxiv23,MV-RLK-SGD-JOTA24}.
The only difference is that the bound \eqref{eq:4321} holds only after
some time $T$.
Clearly this does not affect the \textit{asymptotic} rate of convergence.
Nevertheless, in the interests of completeness, the proof is 
\textit{sketched} here.

The hypotheses on the various constants imply that
\bd
\al_t^2 = O(t^{-2+2 \phi}) , \al_t^2 M_t^2 = O(t^{-2+2 (\phi+\de)}) ,
\al_t B_t = O(t^{-1 + \phi - \ga}) ,
\ed
while $\al_t^2 B_t$ and $\al_t^2 B_t^2$ decay faster than $\al_t B_t$.
Hence both $\{ f_t \}$ and $\{ g_t \}$ are summable if
\bd
-2 + 2 \phi < -1 , -2 + 2(\phi + \de) < -1 , -1 + \phi - \ga < -1 .
\ed
The three inequalities are satisfied if $\phi$ satisfies \eqref{eq:4331a}.
NHext, let us define $\nu$ as in \eqref{eq:4332}, and apply 
Theorem \ref{thm:421}.
This leads to the conclusion that
$J(\w_t) + \nmeusq{\ubold_t} = o(t^{-\la})$ for every $\la \in (0,\nu)$.
In turn this means that, individually, both $J(\w_t)$ and $\nmeusq{\ubold_t}$
are $o(t^{-\la})$ for every $\la \in (0,\nu)$.
Since $\bth_t = \w_t + \ubold_t$, this leads to $J(\bth_t) = o(t^{-\la})$
for every $\la \in (0,\nu)$.
Finally, the (PL) property leads to $\nmeusq{\gJt} = o(t^{-\la})$
for every $\la \in (0,\nu)$.
If we choose the step size sequence to decay very slowly, then the
bound in \eqref{eq:4333} follows readily.
\end{proof}

\subsection{Application to Zero-Order Methods}

In this subsection, we apply Theorem \ref{thm:421} to establish
the convergence of the Stochastc Heavy Ball algorithm when applied to
so-called zero-order (or gradient-free) methods for computing the stochastic
gradient $\h_{t+1}$.
As far back as 1952, a method was introduced in \cite{Kief-Wolf-AOMS52}
for finding a stationary point of a $C^1$ function $J: \R \ap \R$
by approximating the derivative $J'(\theta)$ as
\be\label{eq:451}
J'(\theta) \approx \frac { [J(\theta_t + c_t) + \xi_t^+ ]
- [ J(\theta_t) + \xi_t^- ]}{c_t},
\ee
where $c_t$ is called the ``increment'' at time $t$, and $\xi_t^+ , \xi_t^-$
represent the measurement errors, which are assumed to be i.i.d.\
sequences with zero mean and finite variance.
In \cite{Kief-Wolf-AOMS52} it was observed that, if the above expression
is used as the stochastic gradient $h_{t+1}$, then not only is $h_{t+1}$
a biased estimate of the true derivative $J'(\theta_t)$, but also, the
conditional variance of $h_{t+1}$ is $O(1/c_t^2)$.
In order to make \eqref{eq:451} a better approximation, the increment
$c_t$ is chosen to approach zero as $\tai$.
In turn this causes the conditional variance of $h_{t+1}$ to be an
unbounded function of $t$.
It is shown that, if
\be\label{eq:452}
c_t \ap 0 , \quad \sum_{t=0}^\infty \al_t c_t < \infty, \quad
\sum_{t=0}^\infty (\al_t^2/c_t^2) < \infty , 
\quad \sum_{t=0}^\infty \al_t = \infty .
\ee
then the SGD formulation \eqref{eq:221} converges to a stationary
point of $J(\cdot)$, that is, a solution of $J'(\theta) = 0$.
In \cite{Blum54}, the formulation was extended to functions $J: \R^d \ap \R^d$,
and the convergence of the iterations to a stationary point of $J(\cdot)$
is established under \eqref{eq:452}.
For this reason, the conditions in \eqref{eq:452} are referred to as
the Kiefer-Wolfowitz-Blum conditions, to complement the Robbins-Monro
conditions \eqref{eq:1212}.

Methods such as the above are often called ``zero-order''
or ``gradient-free,'' since they use only function evaluations, and
do not require any gradients to be computed.
As pointed out above, the first such approach is in \cite{Kief-Wolf-AOMS52},
which is shown above as \eqref{eq:451}.
It is for the case $d = 1$, and requires two function
evaluations per iteration.
Subsequently Blum \cite{Blum54} presented an approach for the case $d > 1$,
which requires $d+1$ evaluations per iteration.
When $d$ is large, this approach is clearly impractical.
A significant improvement came in \cite{Spall-TAC92}, in which a method
called ``simultaneous perturbation stochastic approximation'' (SPSA)
was introduced, which requires only \textit{two} function evaluations,
irrespective of the dimension $d$.
However, the proof of convergence of SPSA given in \cite{Spall-TAC92} 
requires many assumptions.
These are simplified in \cite{Chen-Dunc-Dunc-TAC97}.
An ``optimal'' version of SPSA is introduced in \cite{Sadegh-Spall-TAC98},
and is described below.

For each index $t+1$, suppose $\D_{t+1,i}, i \in [d]$ are $d$ different and
pairwise independent
\textbf{Rademacher variables}.\footnote{Recall that Rademacher
random variables assume values in $\bp$ and are independent of each other.}
Moreover, suppose that $\D_{t+1,i} , i \in [d]$ are all independent
(not just conditionally independent) of the $\s$-algebra $\F_t$ for each $t$.
Let $\bD_{t+1} \in \bp^d$ denote the vector of Rademacher variables
at time $t+1$.
Then the search direction $\h_{t+1}$ in \eqref{eq:221} is defined 
componentwise, via
\be\label{eq:453}
h_{t+1,i} = \frac{[ J(\bth_t + c_t \bD_{t+1}) + \xi_{t+1,i}^+ ]
- [ J(\bth_t - c_t \bD_{t+1}) - \xi_{t+1,i}^- ] } {2 c_t \D_{t+1, i}} ,
\ee
where  $\xi_{t+1,1}^+ , \cdots , \xi_{t+1,d}^+$,
$\xi_{t+1,1}^- , \cdots , \xi_{t+1,d}^-$ represent the measurement errors.
A similar idea is used in \cite{Nesterov-FCM17},
except that the bipolar vector $\bD_{t+1}$ is replaced by a random
Gaussian vector $\boldeta_{t+1}$ in $\R^d$.
An excellent survey of this topic can be found in \cite{Li-Xia-Xu-arxiv22},
which discusses other approaches not mentioned here.

The original SPSA envisages only two measurements per iteration, and the
resulting estimate of $\gJ(\bth_t)$ has bias $O(c_t)$ and conditional
variance $O(1/c_t^2)$.
However, it is possible to take more measurements and reduce the bias
of the estimate, while retaining the same bound on the conditional variance.
Specifically, if $k+1$ measurements are taken, then the bias is
$O(c_t^k)$ (which converges to zero more quickly), while the conditional
variance remains as $O(1/c_t^2)$.
See \cite{Pach-Bhat-Pras-arxiv22} and the references therein.

The convergence of the SGD formulation in \eqref{eq:221} is
established in \cite{MV-RLK-SGD-JOTA24}; see specifically
Corollary 6.2.

Against this background, it can be asked whether the Stochastic Heavy Ball
(not SGD) algorithm converges if the stochastic gradient $\h_{t+1}$
is defined as in \eqref{eq:453}.
Note that, even when the measurement errors $\bxi_t^{\pm}$
have zero mean and bounded variance, the stochastic gradient $\h_{t+1}$
defined in \eqref{eq:453} is both biased and has unbounded contitional
variance.
Specifically, if we define $B_t$ and $M_t^2$ as in \eqref{eq:422}
and \eqref{eq:423} respectively, then
\bd
B_t = O(c_t) , M_t = O(1/c_t^2) .
\ed
More generally, if we use the scheme of \cite{Pach-Bhat-Pras-arxiv22}
and use $k+1$ measurements, then
\bd
B_t = O(c_t^k) , M_t = O(1/c_t^2) .
\ed
In either case, previously published papers do not apply to this
situation, especially because of the unbounded variance.
However, Theorem \ref{thm:421} applies to this situation.

\begin{theorem}\label{thm:422}
Consider the Stochastc Heavy Ball Algorithm of \eqref{eq:421},
where the stochastic gradient $\h_{t+1}$ is defined as in \eqref{eq:453},
where $\bxt^\pm$ are zero mean random variables with variance
bounded uniformly with respect to $t$.
Under these conditions, with the same hypotheses as in Theorem
\ref{thm:421}, the conclusions of Theorem \ref{thm:421}
also hold under become the Kiefer-Wolfowitz-Blum conditions of \eqref{eq:452}.
\end{theorem}

The proof is omitted as the stated result is basically a corollary of
Theorem \ref{thm:421}.


\section{A Meta-Theorem on the Convergence of Block Updating}\label{sec:Meta}

Until now we have studied what might be referred to as ``full-coordinate''
updating.
Specifically, in \eqref{eq:130}, \textit{every} coordinate of $\bth_t$
is updated at step $t+1$.
In this section, the focus is on ``block udating,'' wherein, at step $t$,
\textit{some but not necessarily all} components of $\bth_t$ are updated.
Let $S_t \seq [d]$ denote the components of $\bth_t$ that are updated
at step $t$.
Then both the cardinality and the elements of $S_t$ can be random, and
can vary with $t$.
The objective of this section is to prove a ``meta-theorem'' to the
following effect:
Consider thge SGD algorithm of \eqref{eq:130}, and suppose that its
convergence is established using Theorem \ref{thm:23}.
Then the same algorithm, with the same choice of stochastic gradient,
continues to converge under each of three widely used block updating methods.

\subsection{Types of Block Updating Considered}\label{ssec:31}

In this section, we state and prove a meta-theorem for the convergence of 
the SHB algorithm under block updating.
This is achieved by relating the quantities $\z_t$ and $\x_t$
defined in \eqref{eq:422} and \eqref{eq:423} for full coordinate updating
to the corresponding quantities for block updating.
By combining that result with Theorems \ref{thm:421} and
\ref{thm:422} in the present section,
we can directly infer the convergence of SHB with block updating;
no separate proof is required.

Let $\h_{t+1}$ denote the stochastic gradient in \eqref{eq:421}.
The updating method described in \eqref{eq:421} is then the
full coordinate update option.
We refer to it as ``Option 1.''
Now we describe three different options for block updating,
which we call single coordinate, multiple coordinate, and Bernoulli updates.
These are called Options 2, 3 and 4, and are denoted by
$\bphH^{(k)}$ for $k = 2 , 3 , 4$.
These updating schemes
include most if not all of the widely used block updating methods.
In each of these options, while the search direction is random, the
\textit{conditional expectation}
of the update direction is the same as in Option 1.
In symbols, $E_t(\bphH^{(k)}) = E_t(\bphH^{(1)})$ for $k = 2 , 3 , 4$.
In Lemma \ref{lemma:31}, we relate the conditional variance of
$\bphH^{(k)}$ to $\bphH^{(1)}$.
It is worth emphasizing that the conclusions of Lemma \ref{lemma:31}
apply to \textit{arbitrary} stochastic gradients.

Throughout, the symbol $\h_{t+1}$ denotes the search direction
in \eqref{eq:111}.
We now describe Options 1 through 4 for block  updating.

\textbf{Option 1: Full Coordinate Update:}
Let
\be\label{eq:231}
\h^{(1)}_{t+1} = \h_{t+1} .
\ee

\textbf{Option 2: Single Coordinate Update:}
This option is also known as ``coordinate gradient descent'' as defined in
\cite{Wright15} and studied further in \cite{Bach-et-al-aisats19}.
(However, those papers study only the steepest descent method and not 
its variants, such as SHB).
At time $t$, choose an index $\kappa_t \in [d]$ at random with a
uniform probability, and independently of previous choices.
Let $\eb_{\kappa_t}$ denote the elementary unit vector with a $1$
as the $\kappa_t$-th component and zeros elsewhere.
Then define
\be\label{eq:232}
\h^{(2)}_{t+1} = d \eb_{\kappa_t } \circ \h_{t+1} ,
\ee
where $\circ$ denotes the Hadamard, or component-wise, product of two
vectors of equal dimension.
The factor $d$ arises because the likelihood that $\kappa_t$
equaling any one index $i \in [d]$ is $1/d$.

\textbf{Option 3: Multiple Coordinate Update:}
This option is just coordinate update along multiple
coordinates chosen independently at random.
At time $t$, choose $N$ different indices $\kappa_t^n$ from $[d]$
\textit{with replacement}, with each choice being independent of the rest,
and also of past choices.
Moreover, each $\kappa_t^n$ is chosen from $[d]$ with uniform probability.
Then define
\be\label{eq:233}
\h^{(3)}_{t+1} := \frac{d}{N} \sum_{n=1}^N \eb_{\kappa_t^n}
\circ \h_{t+1} .
\ee
Because sampling is with replacement, the average number of times
an index $i \in [d]$ gets selected for updating  is
is $N/d$; to normalize this, the multiplicative factor in \eqref{eq:233} is the
reciprocal of the average.
In this option, $\h^{(3)}_{t+1}$ can have \textit{up to} $N$ nonzero
components.
Because the sampling is \textit{with replacement}, there might
be some duplicated samples.
In such a case, the corresponding component of $\h_{t+1}$ simply
gets counted multiple times in \eqref{eq:233}.

\textbf{Option 4: Bernoulli Update:}
At time $t$, let $\{ B_{t,i}, i \in [d] \}$  be independent Bernoulli
processes with success rate $\rho_t$.
Thus
\be\label{eq:234}
\Pr \{ B_{t,i} = 1 \} = \rho_t , \fa i \in [d] .
\ee
It is permissible for the success probability $\rho_t$ to vary with time.
However, at any one time, all components must have the same success
probability.
Then define
\be\label{eq:235}
\vbold_t := \sum_{i=1}^d \eb_i I_{ \{ B_{t,i} = 1 \} }  \in \bi^d .
\ee
Thus $\vbold_t$ is a random vector, and
$v_{t,i}$ equals $1$ if $ B_{t,i} = 1$, and equals $0$ otherwise.
Now define
\be\label{eq:236}
\h^{(4)}_{t+1} = \frac{1}{\rho_t}  \vbold_t \circ \h_{t+1} .
\ee
Note that, as with the other options,
the factor $1/\rho_t$ is the reciprocal of the likelihood of
a particular $i \in [d]$ being selected for updating.
However, there is no \textit{a priori} upper bound
on the number of nonzero components of $\h^{(4)}_{t+1}$;
the stochastic gradient $\h^{(4)}_{t+1}$ can have up to $d$ nonzero components.
But the \textit{expected} number of nonzero components is $\rho_t d$.

\subsection{A Meta-Theorem on the Convergence of Block Updating}\label{ssec:32}

When the choice of the block update direction involves some
random choices (such as $\kappa_t^n$ or $B_{t+1,i}$), the definition of
the filtration $\{ \F_t \}$ needs to be adjusted.
In the case of Option 2 (coordinate updating), $\F_t$ is the $\s$-algebra
generated by $\kappa_0^t$ in addition to $\bth_0^t$ and $\h_1^t$.
In the case of Option 3, $\kappa_0^t$ is replaced by the collection
$\kappa_{0,i}^t$ for $i \in [N]$.
Finally, in Option 4, $\kappa_0^t$ is replaced by $\vbold_0^t$.

The objectives of Lemma \ref{lemma:31} below are: (i) to show that
all the four search directions have the same conditional expectation,
and (ii) to relate the conditional
variance of Options 2, 3, and 4 to that of Option 1.

\begin{lemma}\label{lemma:31}
As in \eqref{eq:128}, define
\bd
\z_t = E_t(\h_{t+1}) , \bzt = \h_{t+1} - \z_t .
\ed
Then
\be\label{eq:240}
E_t(\h^{(k)}_{t+1} ) = E_t(\bphH^{(1)}) = \z_t , k = 2 , 3 , 4 .
\ee
Moreover,
\be\label{eq:241}
\begin{split}
CV_t( \h^{(2)}_{t+1} ) &= 
(d-1) \nmeusq{\z_t} + d E_t( \nmeusq{\bzt}) , \\
CV_t( \h^{(3)}_{t+1} ) &= 
(d-1) \nmeusq{\z_t} + d E_t( \nmeusq{\bzt}) , \\
CV_t( \h^{(4)}_{t+1} ) &=
\frac{1 -\rho_t}{\rho_t} \nmeusq{\z_t} + \frac{1}{\rho_t} 
E_t( \nmeusq{\bzt} ) .
\end{split}
\ee
\end{lemma}

\begin{proof}
It is obvious that \eqref{eq:240} is satisfied.
Therefore, to compute the conditional variance of $\h^{(k)}_{t+1}$,
it is necessary to compute the residual $\nmeusq{\h^{(k)}_{t+1}  - \z_t}$,
and then take its conditional expectation.

\textbf{Option 2:}
Suppose that $\kappa_t = i$.
Then
\bd
h^{(2)}_{t+1,j} = \left\{ \ba{ll}
d ( z_{t,i} + \zeta_{t+1,i} ) , & \mbox{if } j = i , \\
0, & \mbox{if } j \neq i , \ea \right.
\ed
\bd
h^{(2)}_{t+1,j} - z_{t,j} = \left\{ \ba{ll}
(d-1) z_{t,i} + d \zeta_{t+1,i} , & \mbox{if } j = i , \\
- z_{t,j} , & \mbox{if } j \neq i , \ea \right.
\ed
Therefore, conditioned on the event $\kappa_t = i$, we have that
\bd
\sum_{j=1}^d ( h^{(2)}_{t+1,j} - z_{t,j} )^2 = 
(d-1)^2 z_{t,i}^2 + d^2 \zeta_{t+1,i}^2
+ 2d(d-1)z_{t, i}\zeta_{t+1,i} + \sum_{j\neq i} z_{t,j}^2,
\ed
Now we take the conditional expectation of the above quantity.
For this purpose, we note that (i)
each of the events $\kappa_t = i$ occurs with probability $1/d$,
and (ii) $E_t( z_{t,i} \zeta_{t+1,i} ) = 0$.
Hence
\begin{eqnarray*}
E_t( \nmeusq{ \h^{(2)}_{t+1} - \z_t} ) & = & 
\frac{1}{d}\sum_{i=1}^{d}\left( (d-1)^2 z_{t,i}^2 
+ \sum_{j\neq i} z_{t,j}^2 \right)
+ \frac{1}{d} \sum_{i=1}^d E_t( d^2 \zeta_{t+1,i}^2 ) \\
& = & \frac{(d-1)^2 + (d-1)}{d} \nmeusq{\z_t} + d E_t( \nmeusq{\bzt}) \\
& = & (d-1) \nmeusq{ \z_t } + d E_t(\nmeusq{ \bzeta_{t+1} }).
\end{eqnarray*}
This gives the first equation in \eqref{eq:241}.

\textbf{Option 3:}
Observe that $\h_{t+1}$ is the average of $N$ different quantities
wherein the error terms $\bzeta_{t+1}^n , n \in [N]$ are independent.
Therefore their variances just add up, giving the middle equation in
\eqref{eq:241}.

\textbf{Option 4:}
For notational simplicity, we just use $\rho$ in the place of $\rho_t$.
In this case, each component $h_{t+1,i}$ equals
$ (1/\rho) ( z_{t,i} + \zeta_{t+1,i} )$
with probability $\rho$, and $0$ with probability $1-\rho$.
Thus $h_{t+1,i} - z_{t,i}$ equals $((1/\rho)-1) z_{t,i} + 
(1/\rho) \zeta_{t+1,i}$ with probability $\rho$,
and $- z_{t,i}$ with probability $1-\rho$.
As can be easily verified, the conditional variance is
$((1-\rho)/\rho)z_{t,i}^2 + (1/\rho)E_t(\zeta_{t+1,i}^2))$ for each component.
As the Bernoulli processes for each component are mutually independent,
the variances simply add up.
It follows that
\bd
CV_t( \h^{(4)}_{t+1} ) = \frac{1-\rho}{\rho} \nmeusq{\z_t} 
+ \frac{1}{\rho} E_t( \nmeusq{\bzeta_{t+1} } ) ,
\ed
which is the bottom equation in \eqref{eq:241}.
\end{proof}

With Lemma \ref{lemma:31} in place, we can now state the following meta-theorem
on the convergence of block-uptating applied to the SGD algorithm.

\begin{theorem}\label{thm:441}
Suppose the stochastic gradient $\h_{t+1}$ satisfies the bounds
\eqref{eq:422} and \eqref{eq:423}.
Suppose that in \eqref{eq:130}, the quantity $\h_{t+1}$ is replaced by
$\h_{t+1}^{(k)}$ for $k = 2, 3, 4$.
Further, suppose that when Option $4$ is used, then
\be\label{eq:442}
\inf_t \rho_t =: \bar{\rho} > 0 .
\ee
Then the conclusions of Theorem \ref{thm:23} continue to hold.
\end{theorem}

\begin{proof}
The proof is quite simple, and consists of showing that if $\h_{t+1}$
satisfies \eqref{eq:422} and \eqref{eq:423}, then so do $\h_{t+1}^{(k)}$
for $k = 2, 3, 4$, and then applying Theorem \ref{thm:21}.
In analogy with \eqref{eq:128}, let us define
\bd
\z_t^{(k)} := E_t(\h_{t+1}^{(k)} , \x_t^{(k)} := \z_t^{(k)} - \gJt ,
\bzeta_{t+1}^{(k)} := \h_{t+1}^{(k)} , \mbox{ for } k = 2, 3, 4 .
\ed
Then it follows from \eqref{eq:240} that
\bd
\z_t^{(k)} = \z_t , \x_t^{(k)} = \x_t, \mbox{ for } k = 2, 3, 4 ,
\ed
Now it follows from \eqref{eq:422} that
\bd
\nmeu{\x_t^{(k)}} \leq B_t [ 1 + \nmeu{\gJt} ] ,
\fa \bth_t \in \R^d , \fa t ,k = 2, 3, 4 .
\ed
Next let us prove an analog of \eqref{eq:423} for each $k$.
As a prelude, we can simplify the bounds in \eqref{eq:241} by 
replacing $d-1$ by $d$ and $(1-\rho_t)/\rho_t$ and $1/\rho_t$ by $1/\bar{\rho}$,
where $\bar{\rho}$ is defined in \eqref{eq:442}.
With this substitution, \eqref{eq:241} becomes
\bd
CV_t(\h_{t+1}^{(2)}) \leq d [ \nmeusq{\z_t} + E_t(\nmeusq{\bzt}) ] ,
\ed
\bd
CV_t(\h_{t+1}^{(3)}) \leq d [ \nmeusq{\z_t} + E_t(\nmeusq{\bzt}) ] ,
\ed
\bd
CV_t(\h_{t+1}^{(4)}) \leq \frac{1}{\bar{\rho}}
[ \nmeusq{\z_t} + E_t(\nmeusq{\bzt}) ] .
\ed
With these observations in place, we can simply copy the corresponding
derivation from \cite{MV-RLK-SGD-JOTA24}, specifically the equations
above \cite[Eq.\ (50)]{MV-RLK-SGD-JOTA24}.
Note that the $\mu_t$ in that reference is the present $B_t$ because here
$\mu_t$ denotes the momentum coefficient.
This leads to
\bd
CV_t(\h_{t+1}^{(k)}) \leq d \{ B_t + 2 B_t^2 + [ H(1 + 3 B_t+2 B_t^2 ) + M_t^2 ]
J(\bth_t) \} , \mbox{ for } k = 2 , 3 ,
\ed
\bd
CV_t(\h_{t+1}^{(k)}) \leq  \frac{1}{\bar{\rho}}
\{ B_t + 2 B_t^2 + [ H(1 + 3 B_t + 2 B_t^2 ) + M_t^2 ]
J(\bth_t) \} , k = 4 .
\ed
These bounds are of the form \eqref{eq:213} for suitably defined constants.
Moreover, the analogs of \eqref{eq:222} and \eqref{eq:223} hold with $B_t$
and $M_t$ changed to the new constants, as can be verified easily.
Now the desired convergence follows from Theorem \ref{thm:21}.
\end{proof}

Since the convergence of the SHB algorithm is established by invoking
Theorem \ref{thm:23}, Theorem \ref{thm:441} above implies the
convergence of the SHB algorithm of \eqref{eq:421} under block updating.

\begin{corollary}\label{coro:441}
Suppose the stochastic gradient $\h_{t+1}$ satisfies the bounds
\eqref{eq:422} and \eqref{eq:423}.
Suppose that in \eqref{eq:421}, the quantity $\h_{t+1}$ is replaced by
$\h_{t+1}^{(k)}$ for $k = 2, 3, 4$.
Further, suppose that when Option $4$ is used, then 
$\bar{\rho} > 0$ where $\bar{\rho}$ is defined in \eqref{eq:442}.
Then the conclusions of Theorem \ref{thm:441} continue to hold.
\end{corollary}

\section{Numerical Examples}\label{sec:Num}

In this section, we present several numerical experiments that
illustrate the theory contained in Theorems \ref{thm:423}, \ref{thm:422}
and \ref{thm:441}.
Thus we study the optimization of three different objective functions
(listed below), using stochastic gradients.
The stochastic gradients themselves are of two types: 
(i)  a ``noisy'' gradient, which consists of the true gradient
perturbed by additive Gaussian noise, and (ii) an ``approximate''
gradient using only two function evaluations, as in \eqref{eq:453}.
In the latter case, the stochastic gradient
is biased and its conditional variance grows without bound.
In addition, we study both full-coordinate update as well as block-updating
using Option 4 with Bernoulli sampling.
Though the topic of study here is the stochastic heavy ball (SHB) algorithm,
we also study various other
standard optimization algorithms, such as SGD, SNAG, ADAM, NADAM,
and RMSPROP.\footnote{The author is referred to \cite{Ruder-arxiv16}
for the brief description of these optimization algorithms.}.
In the case of SNAG, we studied two variants: NAG\_F where the step size
$\al_t$ is varied and the momentum term $\mu$ is fixed, and NAG\_S
where the step size is fixed, i.e., $\al_t = al_0$ for all $t$, while
the momentum term $\mu_t$ is scheduled according to the Nesterov's
sequence \cite{Nesterov-Dokl83}.

We used step and increment sequences of the form
\bd
\al_t = \frac{\al_0}{(1 + (t/\tau))^p},
B_t = \frac{b_0}{(1 + (t/\tau))^q}, \mu = 0.9 , \fa t .
\ed
where $\tau$ = 200, $\al_0 = 10^{-6}$, $b_0 = 10^{-4}$, $p = 1$,
and $q = 0.01$.
These parameters are chosen by trial and error, to obtain the best results.
For comparison, we set $\alpha_0$ to be the same
in ADAM, NADAM, and RMSPROP.
For SGD \&  NAG\_F, we used the same $\al_t$ (and $\mu=0$ in case of SGD),
while for NAG\_S, we chose $\al_t = 10^{-6}$ and 
$\mu_t$ to be the Nesterov sequence.
In block updating,
the components to be updated at each time were chosen via
independent Bernoulli processes with a success rate of $\rho$, which was
varied over a range of values.

To evaluate the performance of batch updating, we selected three objective
functions: a strongly convex function
$J_1(\bth_t)$, a non-convex function that satisfies the PL inequality
$J_2(\bth_t)$, and a 2-layer linear neural network loss $J_3(\bth_t)$.
These functions are defined as follows:
\bd
J_1(\bth_t) := \bth_t^\top A \bth_t +
\log \left( \sum_{i=0}^{d-1}e^{\theta_{t, i}} \right),
\ed
\bd
J_2(\bth_t) = \bth_t^\top A \bth_t + 3\sin^2\left(\langle \mathbf{1}, \bth_t \rangle\right),
\ed
\bd
J_3((U_t, V_t)) =
\frac{1}{N}\sum_{i=1}^{N}\nmeusq{\left(M^* - U_tV^\top_t\right)\x_i} +
\lambda \left[\nmm{U_t}_1 + \nmm{V_t}_1\right],
\ed
where $\bth_t$ is a vector of 1 million parameters,
$A$ is a block-diagonal matrix of size ($10^6 \times 10^6$),
consisting of 100 Hilbert matrices
\footnote{The Hilbert matrix is known to be notoriously ill-conditioned,
with eigenvectors that are not aligned with elementary basis vectors
\cite{Fettis-MC67}. This makes it a suitable choice for testing the
robustness of batch updating.}, each of dimension $10^4 \times 10^4$.

In $J_3$, $M^*$ denotes a random target matrix,
$\x_i \sim \N(0, I_{1000 \times 1000})$, and
$U_t \in \R^{1000 \times r}$, $V_t \in \R^{1000 \times r}$ are the learned factors.

Both $J_1$ and $J_2$ can be analytically verified to satisfy the
PL inequality. Empirical results from prior works
\cite{Vidal-AISTATS23, Islamov-et-al-arxiv24} suggest that
$J_3$ also satisfies the PL inequality under certain initialization
and regularity conditions.

Here are the implementation details.
We implemented the
algorithms in \textit{Python} using the \textit{PyTorch} framework.
The experiments were conducted on a workstation equipped with an
Intel Xeon Silver 4114 CPU @ 2.20GHz, 256 GB RAM,
and 4 NVIDIA GeForce RTX 2080 Ti GPUs, each with 12GB of memory.

\begin{figure}[ht]
\centering
\begin{tabular}{cccc}
& {\small  $J_1(\bth_t)$} & {\small $J_2(\bth_t)$} & {\small
 $J_3(\bth_t)$} \\ 
\rotatebox{90}{\small Option 1} &
\includegraphics[width=0.25\linewidth, height=30mm]{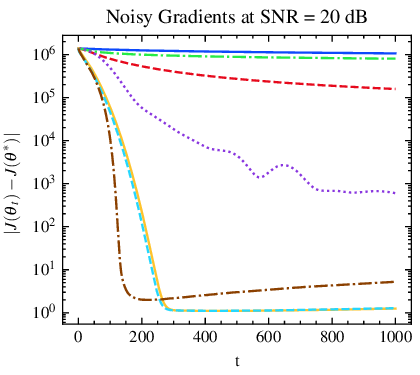} &
\includegraphics[width=0.25\linewidth, height=30mm]{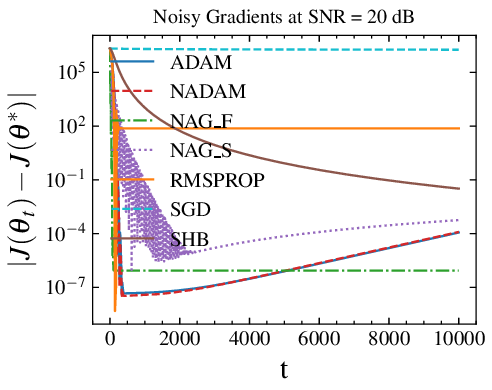} &
\includegraphics[width=0.25\linewidth, height=30mm]{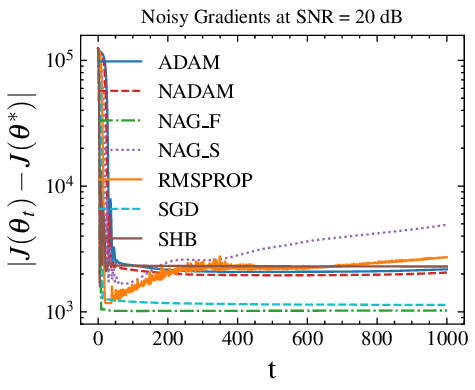} \\
\rotatebox{90}{\small Option 2} &
\includegraphics[width=0.25\linewidth, height=30mm]{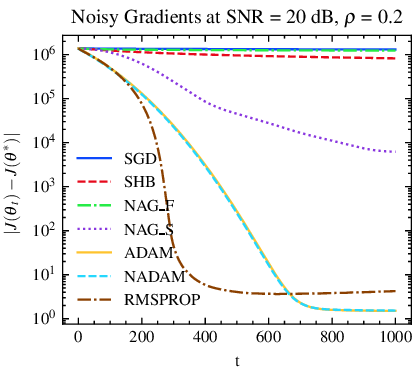} &
\includegraphics[width=0.25\linewidth, height=30mm]{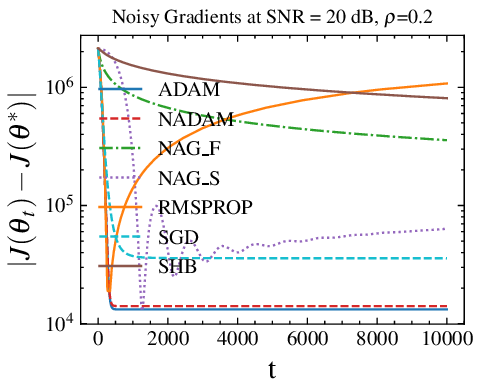} &
\includegraphics[width=0.25\linewidth, height=30mm]{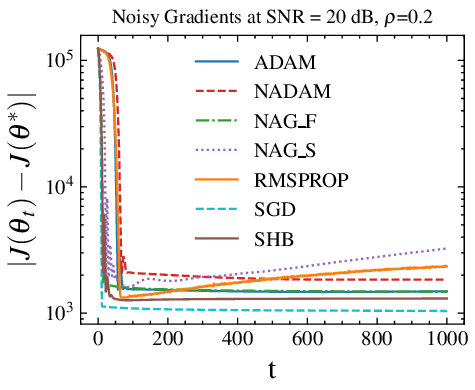} \\
\end{tabular}
\caption{Comparison of various algorithms with noisy gradients (true
gradients corrupted by additive zero mean Gaussian noise)}
\label{fig:exact}
\end{figure}

The results, as shown in Figure \ref{fig:exact}, demonstrate that
when noisy gradients are used, ADAM, NADAM and RMSPROP comfortably
outperform the other four methods.
Within those four, NAG\_S outperforms SHB, which in turn outperforms NAG\_F.
As expected, SGD performs the worst of all.
Moreover, the convergence of ADAM, NADAM, and RMSPROP with Option 4 and
$\rho = 0.2$ (only 20\% of components updated at each iteration) is comparable
to that of full update, after accounting for the reduced updating.

\begin{figure}[ht]
\centering
\begin{tabular}{cccc}
& {\small  $J_1(\bth_t)$} & {\small $J_2(\bth_t)$} & {\small
 $J_3(\bth_t)$} \\ 
\rotatebox{90}{\small Option 1B} &
\includegraphics[width=0.25\linewidth, height=30mm]{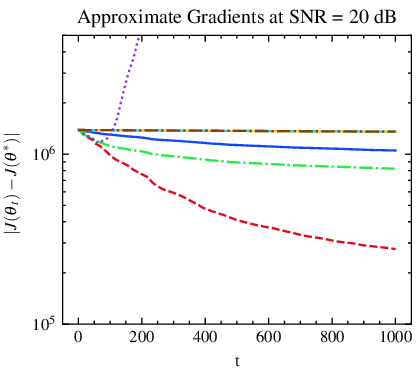} &
\includegraphics[width=0.25\linewidth, height=30mm]{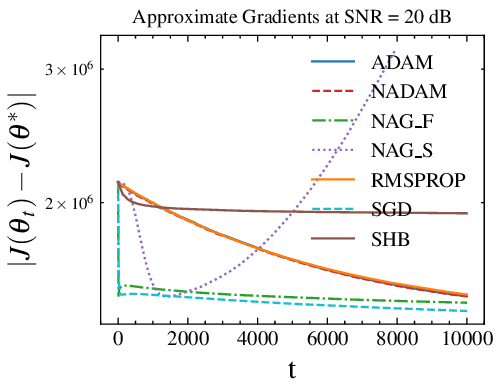} &
\includegraphics[width=0.25\linewidth, height=30mm]{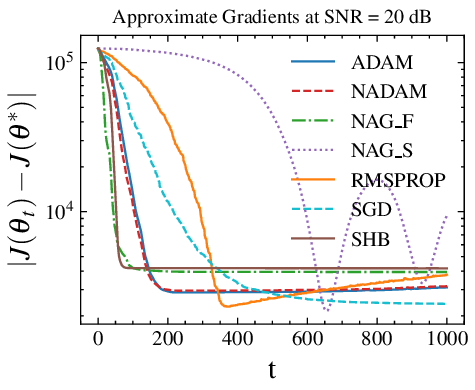} \\
\rotatebox{90}{\small Option 4B} &
\includegraphics[width=0.25\linewidth, height=30mm]{CVX-exact_bcd_final.eps} &
\includegraphics[width=0.25\linewidth, height=30mm]{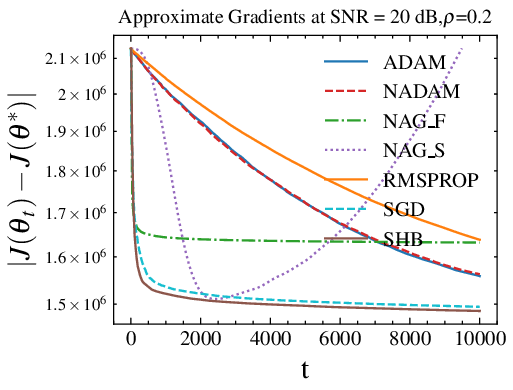} &
\includegraphics[width=0.25\linewidth, height=30mm]{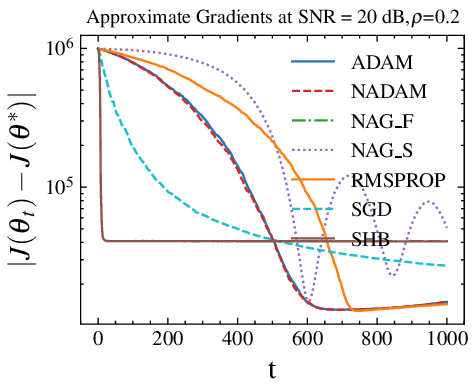} \\
\end{tabular}
\caption{Comparison of various algorithms with approximate gradients
(gradients approximated using \eqref{eq:453})}
\label{fig:approx}
\end{figure}

However, the situation is quite different
when approximate gradients of the form \eqref{eq:453} are used,
As shown in Figure \ref{fig:approx}, NAG\_S diverges almost immediately, 
while ADAM, NADAM, and RMSPROP neither converge nor diverge.
In fact, these three methods perform worse than even SGD.
Among the rest, SGD converges the most slowly, NAG\_F is intermediate,
and SHB performs the best.

Thus, the use of approximate gradients apparently makes it infeasible to use
NAG\_S, ADAM, NADAM, and RMSPROP, whether with full or batch updates.
As pointed out earlier, implementing batch updates with noisy gradient
does not lead to much savings in CPU time because computing only some
components of the gradient vector is almost as CPU-intensive as computing
the entire gradient.
In contrast, with approximate gradients, the amount of computation is
proportional to $\rho$ when Option 4B is used.
The fact that all three methods (SGD, NAG\_F, and SHB) all converge
despite using approximate gradients is very encouraging.

\begin{figure}[ht]
\centering
\begin{tabular}{cc}
\rotatebox{90}{\small $J_1(\bth_1)$} &
\includegraphics[width=0.8\textwidth]{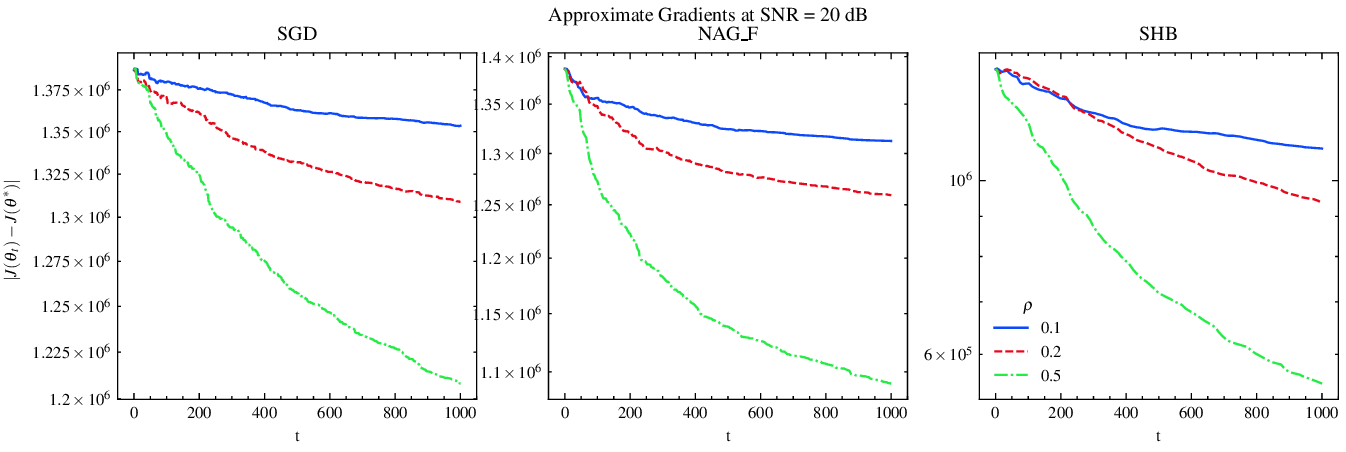}\\
\rotatebox{90}{\small $J_2(\bth_1)$} &
\includegraphics[width=0.8\textwidth]{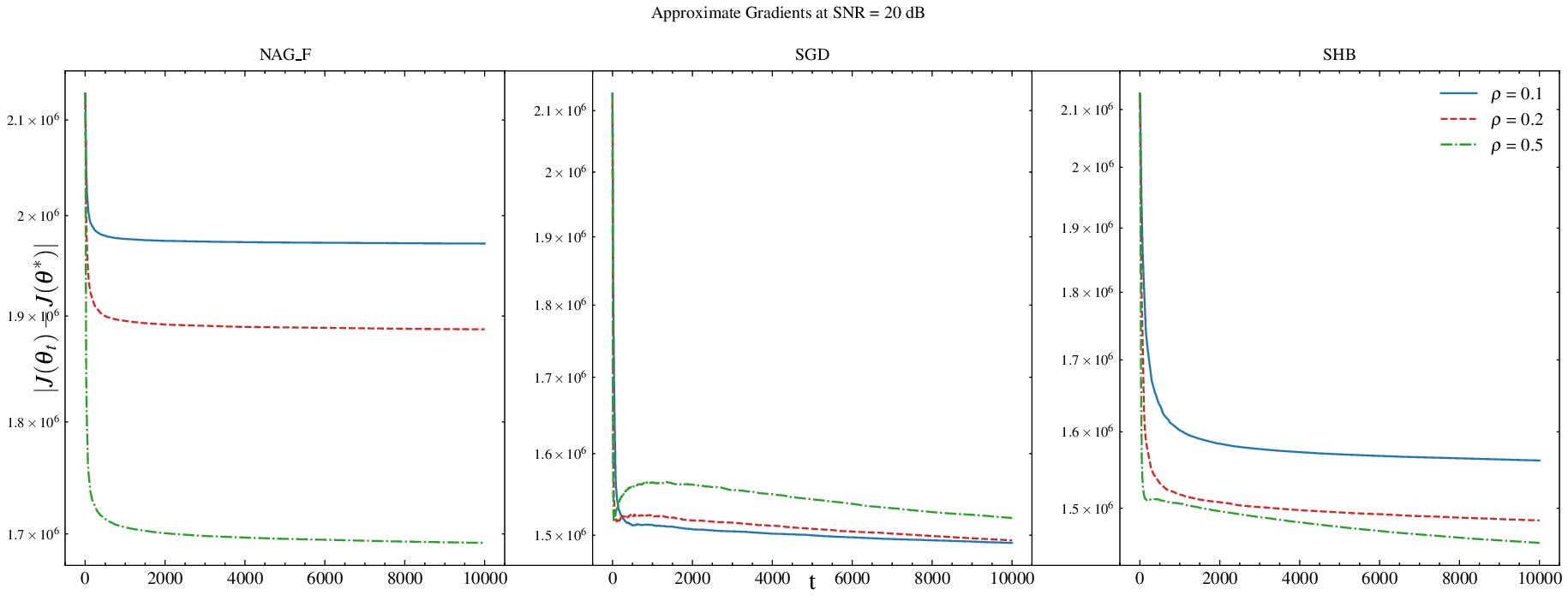}\\
\rotatebox{90}{\small $J_3(\bth_1)$} &
\includegraphics[width=0.8\textwidth]{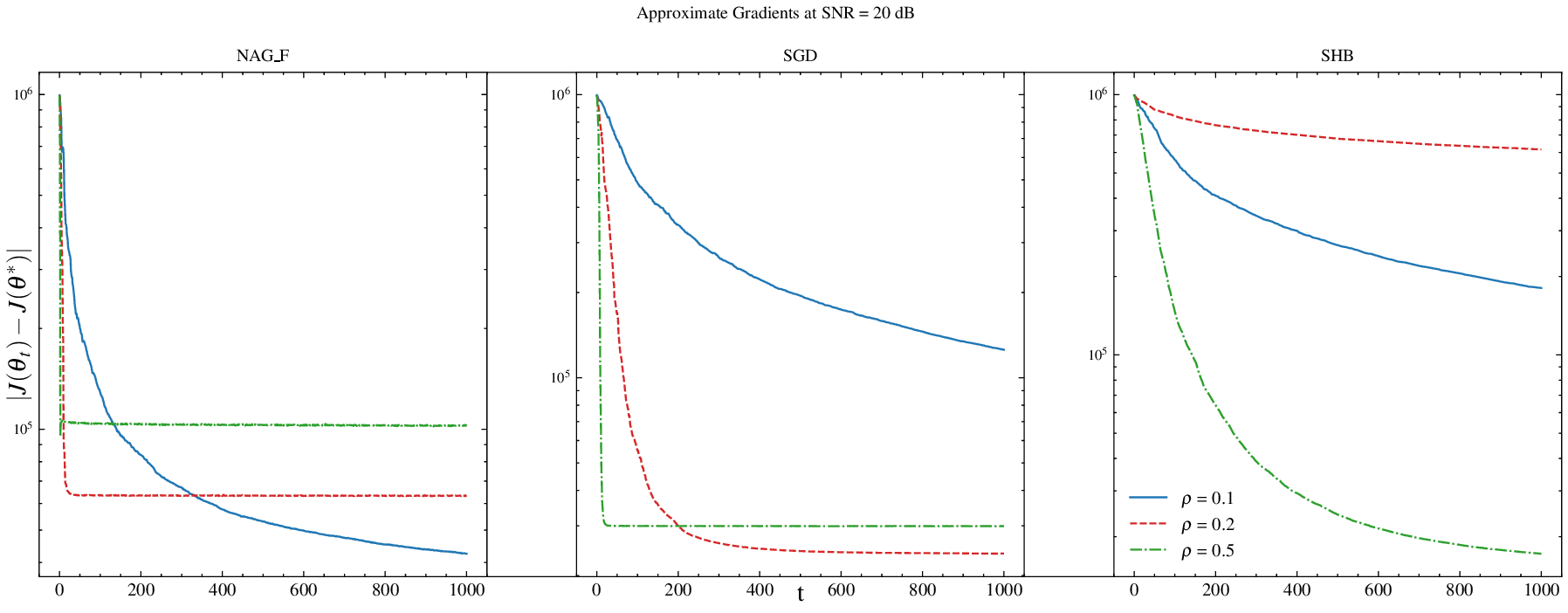}\\
\end{tabular}
\caption{Comparison of various algorithms with noisy gradients
and block updating, with various choices of $\rho$}
\label{fig:varyrho}
\end{figure}

Figure \ref{fig:varyrho} shows that , as expected,
reducing $\rho$ results in slower
convergence because parameters are updated less frequently.
However, the iterations still converge for values of $\rho$ as small as $0.1$,
that is, only 10\% of the components of $\bth_t$ are
updated on average at each iteration.

\begin{figure}[ht]
\centering
\begin{tabular}{cc}
\rotatebox{90}{\small $J_1(\bth_t)$} &
\includegraphics[width=0.8\textwidth]{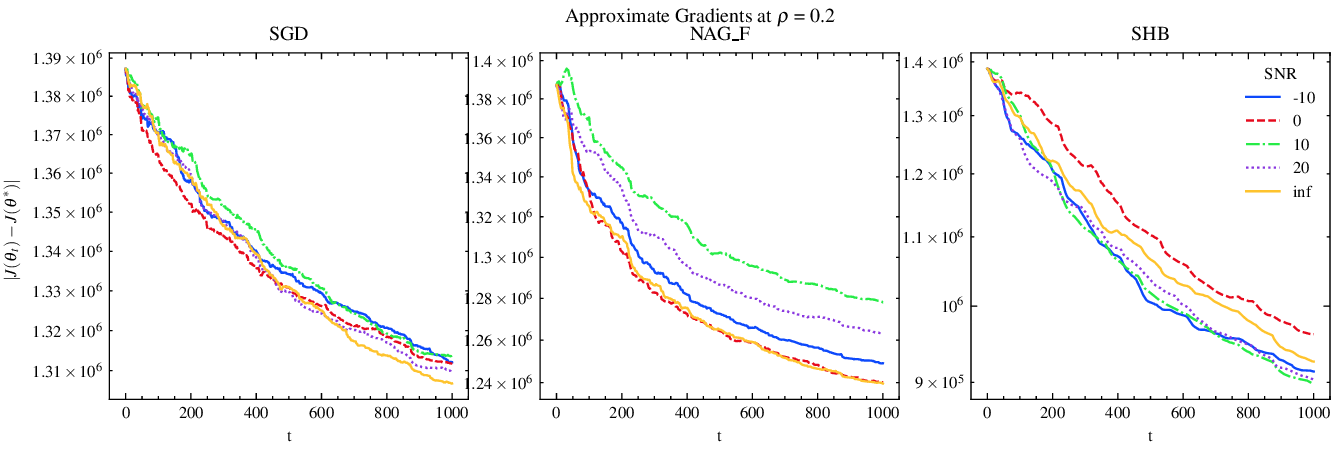}\\
\rotatebox{90}{\small $J_2(\bth_t)$} &
\includegraphics[width=0.8\textwidth]{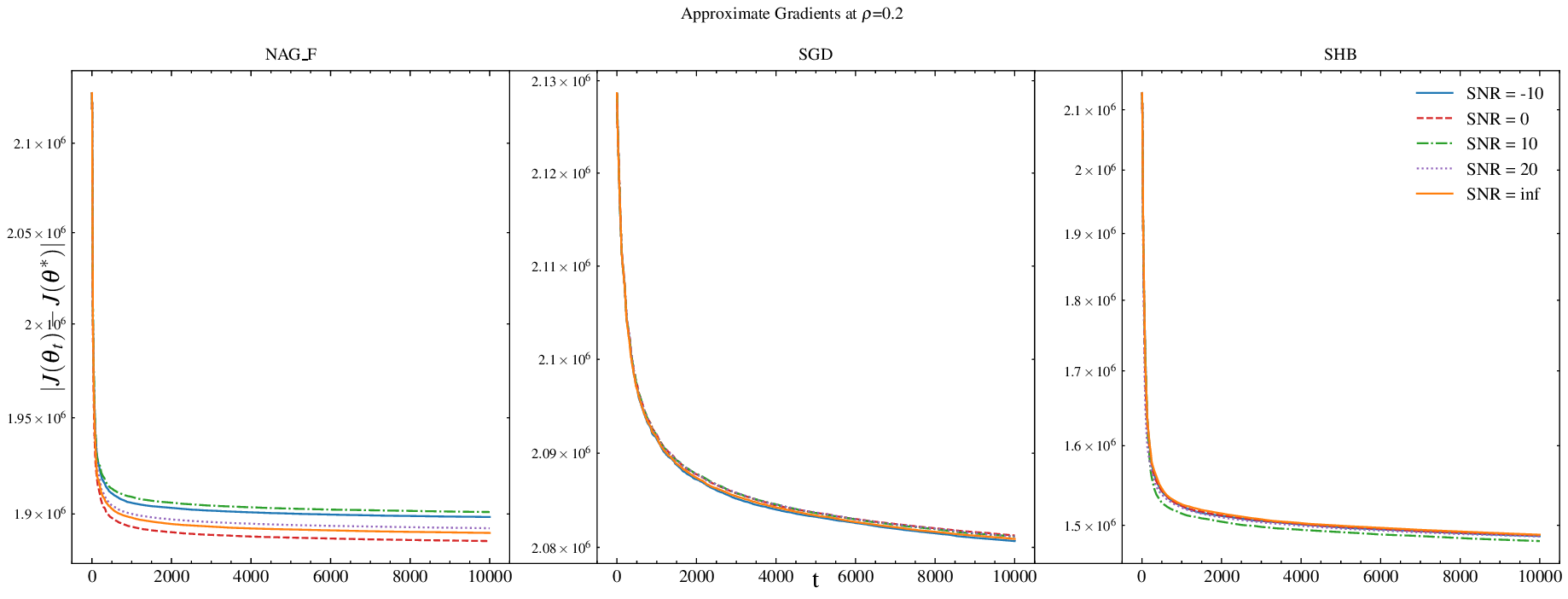}\\
\rotatebox{90}{\small $J_3(\bth_t)$} &
\includegraphics[width=0.8\textwidth]{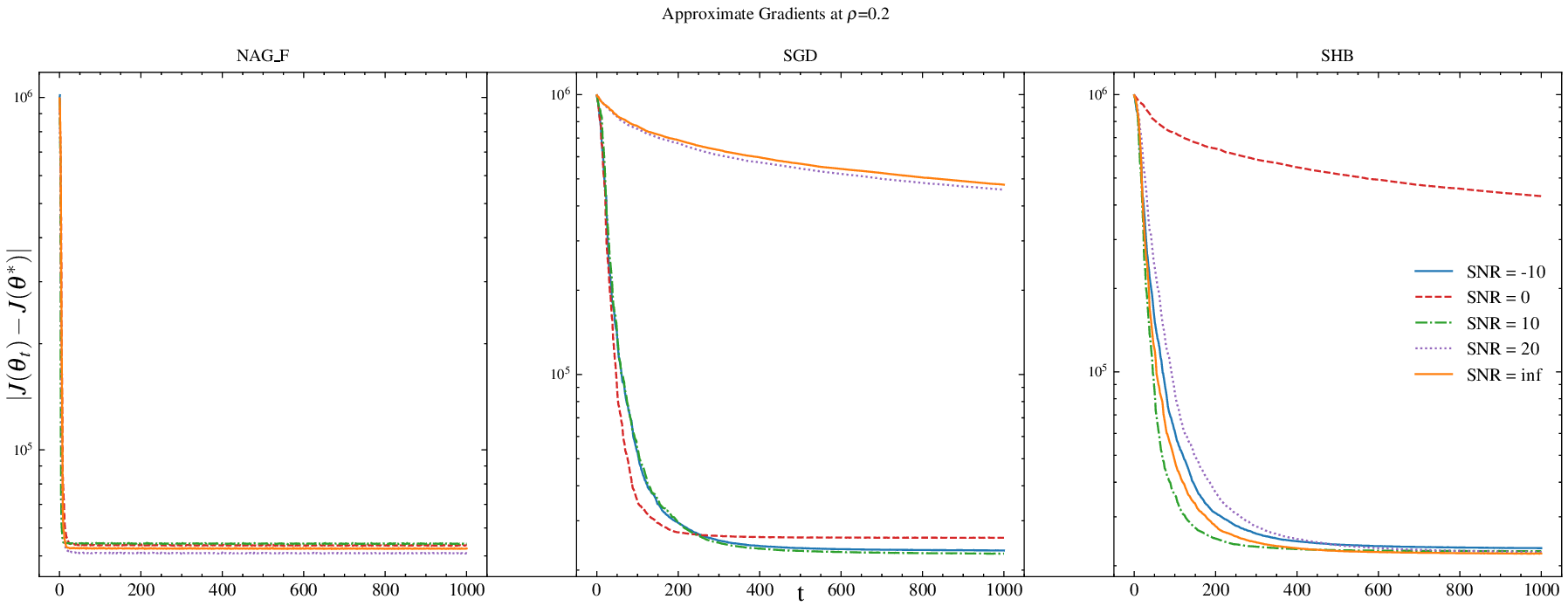}\\
\end{tabular}
\caption{Comparison of various algorithms with approximate gradients
and block updating, with various choices of $\rho$}
\label{fig:varysnr}
\end{figure}

We also tested the robustness
of the batch updating schemes at various noise levels.
The results are shown in Figure \ref{fig:varysnr}. The convergence
rates of these algorithms are comparable at all SNR levels.

\section{Conclusions and Future Work}\label{sec:Conc}

In this paper, we have established the convergence of the stochastic Heavy Ball
(SHB) algorithm under more general conditions than in the current literature.
Specifically,
\bit
\item The stochastic gradient is permitted to be biased, and also, to
have conditional variance that grows over time (or iteration number).
This feature is essential when applying SHB with zeroth-order methods,
which use only two function evaluations to approximate the gradient.
In contrast, all existing papers assume that the stochastic gradient
is unbiased and/or has bounded conditional variance.
\item The step sizes are permitted to be random, which is essential when
applying SHB with block updating.
The sufficient conditions for convergence are stochastic analogs of the
well-known Robbins-Monro conditions.
This is in contrast to existing papers where more restrictive conditions
are imposed on the step size sequence.
\item Our analysis embraces not only convex functions, but also more general
functions that satisfy the PL (Polyak-{\L}ojasiewicz) condition, and KL',
which is slightly weaker than the
(Kurdyka-{\L}ojasiewicz) condition.
\item If the stochastic gradient is unbiased and has bounded variance,
and the objective function satisfies PL), then the iterations of SHB match
the known best rates for convex functions from \cite{Arjevani-et-al-MP23}.
\item We establish the almost-sure convergence of the iterations,
as opposed to convergence in the mean or convergence in probability,
which is the case in much of the literature.
\item Each of the above convergence results continue to hold if 
full-coordinate updating is replaced by any one of three widely-used
updating methods.
\eit
Our current plan is to extend the present results to methods such as
ADAM and RMSPROP, by adapting the methods of \cite{Bara-Bian-SIAM21}.

We have also carried out a series of numerical computations to demonstrate the
following tentative conclusions:
\bit
\item When block updating is applied to noisy gradients, methods such as
ADAM, NADAM, and RMSPROP outperform Stochastic versions of pure gradient
descent, Heavy Ball, and two variants of Nesterov's method.
\item However, when batch updating is applied to approximate gradients,
Nesterov's original method diverges, while ADAM, NADAM, and RMSPROP
barely show any reduction in the objective function.
In fact, they perform worse than the plain steepest descent.
On the other hand, Stochastic Heavy Ball performs the best.
Therefore further theoretical analysis is required to explore
why this is so.
\eit

\section*{Acknowledgements}

The authors thank Mr.\ Lokesh Badisa, final year undergraduate at the
IIT Hyderabad, for assistance with the computations in Section
\ref{sec:Num}.

\section*{Funding and Conflict of Interest}

This work was done when the first author was a final year undergraduate in the
Electrical Engineering Department at IIT Hyderabad.
His work at U.Penn is supported by the Leggett Family Fellowship and the Dean's Fellowship programs.
The research of the second author was supported by the Science and
Engineering Research Board, India.
The authors declare that they do not have any conflict of interest.


\begin{thebibliography}{10}

\bibitem{Apido-et-al-JGO22}
{\sc V.~Apidopoulos, N.~Ginatta, and S.~Villa}, {\em Convergence rates for the
  heavy-ball continuous dynamics for non-convex optimization, under
  polyak–lojasiewicz condition}, Journal of Global Optimization, 84 (2022),
  pp.~563--589.

\bibitem{Arjevani-et-al-MP23}
{\sc Y.~Arjevani, Y.~Carmon, J.~C. Duchi, D.~J. Foster, N.~Srebro, and
  B.~Woodworth}, {\em Lower bounds for non-convex stochastic optimization},
  Mathematical Programming, 199 (2023), pp.~165--214.

\bibitem{Aujol-et-al-SIAMO19}
{\sc J.-F. Aujol, C.~Dossal, and A.~Rondepierre}, {\em Optimal convergence
  rates for nesterov acceleration}, {SIAM} Journal on Optimization, 29 (2019),
  pp.~3131--3153.

\bibitem{Aujol-et-al-MP23}
{\sc J.-F. Aujol, C.~Dossal, and A.~Rondepierre}, {\em Convergence rates of the
  heavy-ball method under the lojasiewicz property}, Mathematical Programming,
  198 (2023), pp.~198--254.

\bibitem{Bach-Moul-NIPS11}
{\sc F.~Bach and E.~Moulines}, {\em Non-asymptotic analysis of stochastic
  approximation algorithms for machine learning}, in Proceedings of the
  Conference on 24th Neural Information Processing Systems (NIPS), 2011,
  pp.~1--9.

\bibitem{Bara-Bian-SIAM21}
{\sc A.~Barakat and P.~Bianchi}, {\em Convergence and dynamical analysis of the
  behavior of the {ADAM} algorithm for nonconvex stochastic optimization},
  {SIAM} Journal on Optimization, 31 (2021), pp.~244--274.

\bibitem{Bengio-et-al-ICASSP13}
{\sc Y.~Bengio, N.~Boulanger-Lewandowski, and R.~Pascanu}, {\em Advances in
  optimizing recurrent networks}, in 2013 IEEE International Conference on
  Acoustics, Speech and Signal Processing, 2013, pp.~8624--8628.

\bibitem{BMP90}
{\sc A.~Benveniste, M.~M\'{e}tivier, and P.~Priouret}, {\em Adaptive Algorithms
  and Stochastic Approximation}, Springer-Verlag, 1990.

\bibitem{Ber-Tsi-SIAM00}
{\sc D.~P. Bertsekas and J.~N. Tsitsiklis}, {\em Global convergence in gradient
  methods with errors}, {SIAM} Journal on Optimization, 10 (2000),
  pp.~627--642.

\bibitem{Blum54}
{\sc J.~R. Blum}, {\em Multivariable stochastic approximation methods}, Annals
  of Mathematical Statistics, 25 (1954), pp.~737--744.

\bibitem{BDLM-TAMS10}
{\sc J.~Bolte, T.~P. Nguyen, J.~Peypouquet, and B.~W. Suter}, {\em
  Characterizations of lojasiewicz inequalities: subgradient flows, talweg,
  convexity}, Transactions of the American Mathematical Society, 362 (2010),
  pp.~3319--3363.

\bibitem{Bottou-et-al-SIAM18}
{\sc L.~Bottou, F.~E. Curtis, and J.~Nocedal}, {\em Optimization methods for
  large-scale machine learning}, SIAM Review, 60 (2018), pp.~223--311.

\bibitem{Chen-Dunc-Dunc-TAC97}
{\sc H.~F. Chen, T.~E. Duncan, and B.~Pasik-Duncan}, {\em A kiefer–wolfowitz
  algorithm with randomized differences}, {IEEE Transactions on Automatic
  Control}, 44 (1999), pp.~442--453.

\bibitem{Chen-et-al-iclr16}
{\sc J.~Chen, R.~Monga, S.~Bengio, and R.~Jozefowicz}, {\em Revisiting
  distributed synchronous sgd}, in International Conference on Learning
  Representations Workshop Track, 2016.

\bibitem{Durrett19}
{\sc R.~Durrett}, {\em Probability: Theory and Examples (5th Edition)},
  Cambridge University Press, 2019.

\bibitem{Fettis-MC67}
{\sc H.~E. Fettis and J.~C. Caslin}, {\em Eigenvalues and eigenvectors of
  hilbert matrices of order 3 through 10}, Mathematics of Computation, 21
  (1967), p.~431.

\bibitem{Fran-Gram21}
{\sc B.~Franci and S.~Grammatico}, {\em Convergence of sequences: A survey}.
\newblock arxiv:2111.11374, 2011.

\bibitem{Gadat-et-al-EJS18}
{\sc S.~Gadat, F.~Panloup, and S.~Saadane}, {\em Stochastic heavy ball},
  Electronic Journal of Statistics, 12 (2018), pp.~461--529.

\bibitem{Ghadimi-et-al-ECC15}
{\sc E.~Ghadimi, H.~R. Feyzmahdavian, and M.~Johansson}, {\em Global
  convergence of the heavy-ball method for convex optimization}, in Proceedings
  of the 2015 European Control Conference, 2015, pp.~311--316.

\bibitem{Ghadimi-Lan-SIAMO13}
{\sc S.~Ghadimi and G.~Lan}, {\em Stochastic first- and zeroth-order methods
  for nonconvex stochastic programming}, {SIAM J. Optimization}, 23 (2013),
  pp.~2341--2368.

\bibitem{Hanson-JMAA81}
{\sc M.~A. Hanson}, {\em On sufficiency of kuhn-tucker conditions}, Journal of
  Mathematical Analysis and Applications, 80 (1981), pp.~545--550.

\bibitem{Hernandez-Spall-ACC14}
{\sc K.~Hernandez and J.~C. Spall}, {\em {Cyclic stochastic optimization with
  noisy function measurements}}, 2014 American Control Conference,  (2014),
  pp.~5204--5209.

\bibitem{Islamov-et-al-arxiv24}
{\sc R.~Islamov, N.~Ajroldi, A.~Orvieto, and A.~Lucchi}, {\em Loss landscape
  characterization of neural networks without over-parametrization}.

\bibitem{MV-RLK-SGD-JOTA24}
{\sc R.~L. Karandikar and M.~Vidyasagar}, {\em Convergence rates for stochastic
  approximation: Biased noise with unbounded variance, and applications},
  Journal of Optimization Theory and Applications, 203 (2024), pp.~2412--2450.

\bibitem{MV-RLK-SGD-arxiv23}
{\sc R.~L. Karandikar and M.~Vidyasagar}, {\em Convergence rates for stochastic
  approximation: Biased noise with unbounded variance, and applications}.
\newblock https://arxiv.org/pdf/2312.02828v3.pdf, May 2024.

\bibitem{Karimi-et-al16}
{\sc H.~Karimi, J.~Nutini, and M.~Schmidt}, {\em Linear convergence of gradient
  and proximal-gradient methods under the polyak- lojasiewicz condition},
  Lecture Notes in Computer Science, 9851 (2016), pp.~795--811.

\bibitem{Khaled-Rich-arxiv20}
{\sc A.~Khaled and P.~Richt\'{a}rik}, {\em {Better Theory for SGD in the
  Nonconvex World}}.
\newblock arXiv:2002.03329, February 2020.

\bibitem{Kief-Wolf-AOMS52}
{\sc J.~Kiefer and J.~Wolfowitz}, {\em Stochastic estimation of the maximum of
  a regression function}, Annals of Mathematical Statistics, 23 (1952),
  pp.~462--466.

\bibitem{Kurdyka98}
{\sc K.~Kurdyka}, {\em On gradients of functions definable in o-minimal
  structures}, Annales de L'Institut Fourier, 48 (1998), pp.~769--783.

\bibitem{Li-Xia-Xu-arxiv22}
{\sc S.~Li, Y.~Xia, and Z.~Xu}, {\em Simultaneous perturbation stochastic
  approximation: towards one-measurement per iteration}.
\newblock arxiv:2203.03075, March 2022.

\bibitem{Liu-et-al-NC23}
{\sc J.~Liu, D.~Xu, Y.~Lu, J.~Kong, and D.~P. Mandic}, {\em Last-iterate
  convergence analysis of stochastic momentum methods for neural networks},
  Neurocomputing, 527 (2023), pp.~27--35.

\bibitem{Liu-Yuan-arxiv22}
{\sc J.~Liu and Y.~Yuan}, {\em On almost sure convergence rates of stochastic
  gradient methods}, in Proceedings of Thirty Fifth Conference on Learning
  Theory, P.-L. Loh and M.~Raginsky, eds., vol.~178 of Proceedings of Machine
  Learning Research, PMLR, 02--05 Jul 2022, pp.~2963--2983.

\bibitem{Loja63}
{\sc S.~Lojasiewicz}, {\em Une propri\'{e}t\'{e} topologique des sous-ensembles
  analytiques r\'{e}els}, Editions du centre National de la Recherche
  Scientifique, 1963, pp.~87--89.

\bibitem{Lu-Xiao-arxiv13}
{\sc Z.~Lu and L.~Xiao}, {\em On the complexity analysis of randomized
  block-coordinate descent methods}, Mathematical Programming, 152 (2014),
  pp.~615--642.

\bibitem{Mertik-et-al-NeurIPS20}
{\sc P.~Mertikopoulos, N.~Hallak, A.~Kavis, and V.~Cevher}, {\em On the almost
  sure convergence of stochastic gradient descent in non-convex problems}, in
  Advances in Neural Information Processing Systems, H.~Larochelle, M.~Ranzato,
  R.~Hadsell, M.~Balcan, and H.~Lin, eds., vol.~33, 2020, pp.~1117--1128.

\bibitem{Nemirovski-et-al-SIAMO09}
{\sc A.~Nemirovski, A.~Juditsky, G.~Lan, and A.~Shapiro}, {\em Robust
  stochastic approximation approach to stochastic programming}, {SIAM Journal
  on optimization}, 19 (2009), pp.~1574--1609.

\bibitem{Nesterov-Dokl83}
{\sc Y.~Nesterov}, {\em A method for unconstrained convex minimization problem
  with the rate of convergence $o(1/k^2)$n (in russian)}, Soviet Mathematics
  Doklady, 269 (1983), pp.~543--547.

\bibitem{Nesterov-SICOPT12}
{\sc Y.~Nesterov}, {\em Efficiency of coordinate descent methods on huge-scale
  optimization problems}, SIAM Journal on Optimization, 22 (2012),
  pp.~341--362.

\bibitem{Nesterov18}
{\sc Y.~Nesterov}, {\em Lectures on Convex Optimization (Second Edition)},
  Springer Nature, 2018.

\bibitem{Nesterov-FCM17}
{\sc Y.~Nesterov and V.~Spokoiny}, {\em {Random Gradient-Free Minimization of
  Convex Functions}}, Foundations of Computational Mathematics, 17 (2017),
  pp.~527--566.

\bibitem{Pach-Bhat-Pras-arxiv22}
{\sc S.~Pachal, S.~Bhatnagar, and L.~A. Prashanth}, {\em Generalized
  simultaneous perturbation-based gradient search with reduced estimator bias}.
\newblock arxiv:2212.10477v1, December 2022.

\bibitem{Polyak-CMMP64}
{\sc B.~Polyak}, {\em Some methods of speeding up the convergence of iteration
  methods}, USSR Computational Mathematics and Mathematical Physics, 4 (1964),
  pp.~1--17.

\bibitem{Polyak-UCMMP63}
{\sc B.~T. Polyak}, {\em Gradient methods for the minimisation of functionals},
  {USSR Computational Mathematics and Mathematical Physics}, 3 (1963),
  pp.~864--878.

\bibitem{Pol-Jud-SICOPT92}
{\sc B.~T. Polyak and A.~B. Juditsky}, {\em Acceleration of stochastic
  approximation by averaging}, {SIAM Journal of Control and Optimization}, 30
  (1992), pp.~838--855.

\bibitem{Niu-et-al-arxiv11}
{\sc B.~Recht, C.~Re, S.~Wright, and F.~Niu}, {\em Hogwild!: A lock-free
  approach to parallelizing stochastic gradient descent}, in Advances in Neural
  Information Processing Systems, J.~Shawe-Taylor, R.~Zemel, P.~Bartlett,
  F.~Pereira, and K.~Weinberger, eds., vol.~24, Curran Associates, Inc., 2011.

\bibitem{Richtarik-Takac-MP12}
{\sc P.~Richt{\'{a}}rik and M.~Tak{\'{a}}{\v{c}}}, {\em Iteration complexity of
  randomized block-coordinate descent methods for minimizing a composite
  function}, Mathematical Programming, 144 (2012), pp.~1--38.

\bibitem{Richtarik-Takac-MP15}
{\sc P.~Richt{\'{a}}rik and M.~Tak{\'{a}}{\v{c}}}, {\em Parallel coordinate
  descent methods for big data optimization}, Mathematical Programming, 156
  (2015), pp.~433--484.

\bibitem{Robbins-Monro51}
{\sc H.~Robbins and S.~Monro}, {\em A stochastic approximation method}, Annals
  of Mathematical Statistics, 22 (1951), pp.~400--407.

\bibitem{Robb-Sieg71}
{\sc H.~Robbins and D.~Siegmund}, {\em A convergence theorem for non negative
  almost supermartingales and some applications}, Elsevier, 1971, pp.~233--257.

\bibitem{Ruder-arxiv16}
{\sc S.~Ruder}, {\em An overview of gradient descent optimization algorithms},
  2016.

\bibitem{Sadegh-Spall-TAC98}
{\sc P.~Sadegh and J.~C. Spall}, {\em Optimal random perturbations for
  stochastic approximation using a simultaneous perturbation gradient
  approximation}, {IEEE transactions on automatic control}, 43 (1998),
  pp.~1480--1484.

\bibitem{Sebbouh-et-al-CoLT21}
{\sc O.~Sebbouh, R.~M. Gower, and A.~Defazio}, {\em Almost sure convergence
  rates for stochastic gradient descent and stochastic heavy ball}, in
  Proceedings of Thirty Fourth Conference on Learning Theory, {PMLR}, vol.~134,
  2021, pp.~3935--3971.

\bibitem{Spall-TAC92}
{\sc J.~Spall}, {\em {Multivariate stochastic approximation using a
  simultaneous perturbation gradient approximation}}, IEEE Transactions on
  Automatic Control, 37 (1992), pp.~332--341.

\bibitem{Spall-TAES98}
{\sc J.~C. Spall}, {\em Implementation of the simultaneous perturbation
  algorithm for stochastic optimization}, {IEEE transactions on aerospace and
  electronic systems}, 34 (1998), pp.~817--823.

\bibitem{Su-Boyd-Candes16}
{\sc W.~Su, S.~Boyd, and E.~J. Cand\`{e}s}, {\em A differential equation for
  modeling nesterov’s accelerated gradient method: Theory and insights},
  Journal of Machine Learning Research, 17 (2016), pp.~1--43.

\bibitem{Sutskever-et-al-ICML13}
{\sc I.~Sutskever, J.~Martens, G.~Dahl, and G.~Hinton}, {\em On the importance
  of initialization and momentum in deep learning}, in Proceedings of the 30th
  international conference on machine learning (ICML-13), 2013, pp.~1139--1147.

\bibitem{Bach-et-al-aisats19}
{\sc S.~Vaswani, F.~Bach, and M.~Schmidt}, {\em {Fast and Faster Convergence of
  SGD for Over-Parameterized Models (and an Accelerated Perceptron)}}, in
  {Proceedings of the 22nd International Conference on Artificial Intelligence
  and Statistics (AISTATS)}, 2019, pp.~1--10.

\bibitem{Williams91}
{\sc D.~Williams}, {\em Probability with Martingagles}, Cambridge University
  Press, 1991.

\bibitem{Wright15}
{\sc S.~J. Wright}, {\em Coordinate descent algorithms}, Mathematical
  Programming, 151 (2015), pp.~3--34.

\bibitem{Xie-et-al-ICML20}
{\sc C.~Xie, S.~Koyejo, and I.~Gupta}, {\em Zeno++: Robust fully asynchronous
  {SGD}}, in Proceedings of the 37th International Conference on Machine
  Learning, H.~D. III and A.~Singh, eds., vol.~119 of Proceedings of Machine
  Learning Research, PMLR, 13--18 Jul 2020, pp.~10495--10503.

\bibitem{Xu-Yin-SICOPT15}
{\sc Y.~Xu and W.~Yin}, {\em Block stochastic gradient iteration for convex and
  nonconvex optimization}, SIAM Journal on Optimization, 25 (2015),
  pp.~1686--1716.

\bibitem{Vidal-AISTATS23}
{\sc Z.~Xu, H.~Min, S.~Tarmoun, E.~Mallada, and R.~Vidal}, {\em Linear
  convergence of gradient descent for finite width over-parametrized linear
  networks with general initialization}, in Proceedings of The 26th
  International Conference on Artificial Intelligence and Statistics, 2023,
  pp.~2262--2284.

\end{thebibliography}
\end{document}